%NOTE: This document requires AMSTeX 2.0, the 
%accompanying version of the amsppt style macros,
%and the font package amsfonts 2.0. If you try 
%to typeset it using an older version of TeX, it 
%probably won't work.  
%
%%%%%%%%%%%%%%%%%%%  BEGINNING OF cpdilation.tex  %%%%%%%%%%%
%
\input amstex
\documentstyle{amsppt}
\loadbold
\def\cstar{$C^*$-algebra}
\def\esg{$E_0$-semigroup}

\def\<{\left<}										%for inner products
\def\>{\right>}

\magnification=\magstep 1

\topmatter
\title On the index and dilations of completely positive semigroups
\endtitle

\author William Arveson
\endauthor

\affil Department of Mathematics\\
University of California\\Berkeley CA 94720, USA
\endaffil

\date 28 August 1996
\enddate
\thanks This research was supported by
NSF grant DMS95-00291
\endthanks
\keywords index, quantum dynamical semigroups, 
\esg s, von Neumann algebras, completely positive maps
\endkeywords
\subjclass
Primary 46L40; Secondary 81E05
\endsubjclass
\abstract
It is known that 
every semigroup of normal completely positive maps
$P = \{P_t: t\geq 0\}$ of $\Cal B(H)$, satisfying
$P_t(\bold 1) = \bold 1$ for every $t\geq 0$, 
has a minimal 
dilation to an \esg\ acting on $\Cal B(K)$ for 
some Hilbert space $K\supseteq H$.  
The minimal dilation of $P$ is 
unique up to conjugacy.  In a previous paper a 
numerical index was introduced for semigroups 
of completely positive maps and it was shown that
the index of $P$ agrees with the index of its 
minimal dilation to an \esg.  However, no examples
were discussed, and no computations were made.  

In this paper we calculate the index 
of a unital completely positive semigroup whose
generator is a {\it bounded} operator 
$$
L:\Cal B(H)\to\Cal B(H)  
$$
in terms of natrual structures associated with 
the generator.  This includes all unital CP semigroups
acting on matrix algebras.  We also show that the 
minimal dilation of the semigroup $P=\{\exp{tL}: t\geq 0\}$ 
to an \esg\ is is cocycle 
conjugate to a $CAR/CCR$ flow.  
\endabstract

\endtopmatter
%
%\vfill\eject
%Replace \pagebreak below with the line above
%to fill the lower part of the title page with
%space, rather than stretching it.
%\pagebreak

\document

\subheading{Introduction}

In \cite{4}, a numerical index is introduced for semigroups 
$P=\{P_t: t\geq 0\}$ of normal 
completely positive maps of $\Cal B(H)$.  
In the case where $P_t(\bold 1)=\bold 1$ for every $t$, 
a recent theorem of B. V. R. Bhat asserts that $P$ can 
be ``dilated" to an \esg\ \cite{5,6}; and it was shown in 
\cite{4} that the index of $P$ agrees with the index of 
its {\it minimal} \esg\ dilation.  
However, no examples were discussed there.  
In particular, the results of 
\cite{4} give no information 
about which \esg s can occur as the minimal dilations of 
unital completely positive semigroups acting on matrix 
algebras.  

In this paper we consider the more general case of 
completely positive semigroups having
bounded generator.  We calculate the 
index of such semigroups in terms of basic 
structures associated with their generators (Theorem 2.3, 
Corollary 2.17) and in the case where the semigroup 
preserves the unit we show that their minimal dilations must
be cocycle conjugate to a $CAR/CCR$ flow (Corollary 4.21).  
The extent to which the index calculations of section 2 can be 
extended to semigroups with unbounded generators remains
unclear at present; but certainly the description of 
their minimal dilations (e.g., Corollary 4.21) becomes
 false without strong 
hypotheses on the generator.  

It is appropriate to point out that, using a completely 
different method, Powers \cite{11} has independently shown 
that every unital completely 
positive semigroup acting on a matrix algebra dilates to 
a completely spatial \esg\ and he calculates 
the index of the minimal dilation in that case.

\subheading{1. Bounded generators, symbols, and metric operator spaces}

We are concerned with the structure of various linear 
mappings on the von Neumann algebra 
$M = \Cal B(H)$, $H$ being a separable Hilbert space.  
$\Cal L(M)$ will denote the space of all bounded 
linear maps $L: M\to M$.  
The purpose of this section is to discuss the 
relationship of metric operator spaces \cite{4} 
to the generators of completely positive semigroups 
and their symbols.  Our methods in \S\S 2--3, even 
the statement of the key Theorem 2.3, will involve
metric operator spaces in an essential way.

We briefly recall the definition
of the symbol of a linear map $L\in \Cal L(M)$.  
Consider the bilinear mapping 
defined on $M\times M$ by 
$$
L(xy) - xL(y) - L(x)y + xL(\bold 1)y. \tag{1.1} 
$$
It is useful to regard this as a homomorphism 
of the bimodule $\Omega^2$ of all noncommutative 2-forms 
into $M$, and that homomorphism of $M$-modules is 
the symbol of $L$.  More explicitly, 
$\Omega^1$ is defined as the submodule of 
the symmetric bimodule $M\otimes M$ (with operations 
$a(x\otimes y)b = ax\otimes yb$, 
$(x\otimes y)^* = y^*\otimes x^*$) generated by the 
range of the derivation $d: M\to M\otimes M$, 
$$
dx = x\otimes \bold 1 - \bold 1\otimes x.  
$$
We have $(dx)^* = -d(x^*)$ and every element of $\Omega^1$ 
is a finite sum of the form $adx_1 + \dots + adx_n$, 
$a_k, x_k\in M$.  $\Omega^2$ is defined by 
$$
\Omega^2 = \Omega^1\otimes_M\Omega^1.  
$$
Every element of $\Omega^2$ is a sum of 
the form $a_1dx_1\,dy_1 + \dots + a_ndx_n\,dy_n$, 
and we have a natural multiplication
$$
\omega_1, \omega_2\in \Omega^1 \to 
\omega_1 \omega_2\in \Omega^2
$$
which satisfies the associative law 
$\omega_1(a\omega_2) = (\omega_1 a)\omega_2$ with
respect to operators $a\in M$.  

Given $L\in \Cal L(M)$ there is a unique 
$\sigma_L\in\hom(\Omega^2,M)$ satisfying 
$$
\sigma_L(dx\,dy) = L(xy)-xL(y) -L(x)y + xL(\bold 1)y
$$
(see \cite{3}).  $\sigma_L$ is called the {\bf symbol}
of $L$, and it has the following basic
properties:  

\roster
\item"{(1.2)}"
$\sigma_L = 0$ iff $L$ has the form $L(x) = ax + xb$ for
fixed elements $a,b\in M$. 
\item"{(1.3)}"
 If $\sigma_L=0$ and $L$ satisfies $L(x^*) = L(x)^*$ for every
$x$, then there is an element $a\in M$ such that 
$L(x) = ax + xa^*$.  
\item"{(1.4)}" 
For every $x,y\in M$ we have 
$$
\|\sigma_L(dx\,dy)\| \leq 4\|L\|\,\|x\|\,\|y\|.  
$$
\endroster
\remark{Remark 1.5}
Property (1.2) follows from the fact that if 
$\sigma_L = 0$ then the linear map 
$L_0(x) = L(x) - xL(\bold 1)$ is a derivation of $M$, 
and hence has the form $L_0(x) = ax - xa$.  Property 
(1.3) follows from (1.2) after an elementary 
argument which uses the fact that if $c$ is any operator 
satisfying $cx + xc^* = 0$ for every operator $x$
then $c$ must have the form 
$c=\sqrt{-1}\lambda\bold 1$ where $\lambda\in \Bbb R$
(see Lemma 1.19 below).
\endremark

The following proposition summarizes some known results
that illustrate how properties of the symbol characterize
the generators of semigroups of completely positive
maps.  

\proclaim{Proposition 1.6}
Suppose that $A$ is a unital \cstar\ and 
$L\in \Cal L(A)$ is a bounded operator satisfying 
$L(x)^* = L(x^*)$, $x\in A$.  Let $\{P_t: t\geq 0\}$
be the semigroup of linear operators on $A$ obtained by 
exponentiation: $P_t = \exp(tL)$.  Then the following 
are equivalent.  
\roster
\item"{1.6.1}" Each map $P_t$ is completely positive.   
\item"{1.6.2}" If $a_1, x_1, \dots, a_n, x_n\in A$ satisfy 
$x_1a_1 + \dots x_na_n=0$ then we have
$$
\sum_{k,j=1}^n a_j^*L(x_j^*x_k)a_k \geq 0.  
$$
\item"{1.6.3}"$\sigma_L(\omega^*\omega)\leq 0$ 
for all $\omega_1,\omega_2\in \Omega^1$.  
\endroster
\endproclaim

\demo{proof}
The equivalence of (1.6.1) and (1.6.2) is essentially a result
of Evans and Lewis \cite{8}.  The equivalence of (1.6.2) 
and (1.6.3) is discussed in \cite{3}.  
\qed
\enddemo

It is a straightforward consequence 
of Stinespring's theorem that every
normal completely positive linear map $P\in \Cal L(M)$ can 
be expressed in the form 
$$
P(x) = \sum_kv_kxv_k^* \tag{1.7}
$$
where $\{v_1, v_2,\dots\}$ is a (finite or infinite) sequence
of operators in $M$.  Since certain facts relating to this
representation are fundamental to our approach, 
we offer the following comments.  
Let $(\pi, V)$ be a pair consisting
of a representation $\pi$ of $M=\Cal B(H)$ on some other Hilbert
space space $K$ and an operator $V\in \Cal B(H,K)$ satisfying 
$$
P(x) = V^*\pi(x)V.  \tag{1.8} 
$$
By cutting down to a subspace of $K$ if necessary we can assume
that $K$ is spanned by $\{\pi(x)\xi: x\in M, \xi\in H\}$, and 
in this case the normality of $P$ implies that $\pi$ is a normal
representation, necessarily nondegenerate.  
Thus by replacing $(\pi, V)$ by an equivalent pair we may assume
that $K = H^n$ is a countable 
direct sum of copies of $H$ and $\pi$ has 
the form $\pi(x) = x\oplus x\oplus \dots $.  It follows that
there is a sequence $v_1, v_2, \dots\in \Cal B(H)$ such that
$$
V\xi = (v_1\xi, v_2\xi, \dots)
$$
and the representation (1.7) follows by taking $u_k=v_k^*$.  

There is a natural operator space $\Cal E_P$ associated
with $P$, which can be defined in concrete terms as follows.  
Notice that if $\lambda = (\lambda_1, \lambda_2,\dots)$ is 
any sequence in $\ell^2$ then the operator sum 
$\sum_k\lambda_ku_k$ is convergent in the strong operator 
topology (this sum represents the composition of $V^*$ with
the operator 
$\xi\in H\mapsto (\lambda_1\xi, \lambda_2\xi, \dots)\in H^n$), 
and because of the minimality of $(\pi,V)$ we have 
$$
\lambda_1u_1 + \lambda_2u_2\dots = 0\implies 
\lambda_1 = \lambda_2 = \dots = 0,  \tag{1.9}
$$
for every $\lambda\in \ell^2$.  We define
$$
\Cal E_P = 
\{\lambda_1u_1 + \lambda_2 u_2 + \dots: \lambda\in \ell^2\}.  
\tag{1.10}
$$
$\Cal E_P$ is not necessarily closed in the operator norm
when it is infinite dimensional, but in all cases it 
is a Hilbert space with respect to the inner product 
defined on it by declaring $\{u_1, u_2, \dots \}$ to 
be an orthonormal basis.  $\Cal E_P$ has the following 
properties.  
\roster
\item"{(1.11.1)}"
An operator $a$ belongs to $\Cal E_P$ if and only if 
there is a constant $c\geq 0$ such that the mapping
$$
x\in M \mapsto c P(x) - axa^*
$$
is completely positive.  In this case, $\<a,a\>_\Cal E$ 
is the smallest such constant $c$.   
\item"{(1.11.2)}"
If $w_1, w_2, \dots$ is any orthonormal basis for $\Cal E_P$
then the sum $\sum_kw_kw_k^*$ converges strongly, and in fact 
$$
P(x) = \sum_k w_kxw_k^*, \qquad x\in M.  
$$
\item"{(1.11.3)}"
If $z_1, z_2, \dots$ is any finite or infinite sequence of 
operators in $M$ such that the series $\sum_kz_kz_k^*$ converges
strongly, and which represents $P$ in the 
sense that 
$$
P(x) = \sum_kz_kxz_k^*, \qquad x\in M,
$$
then $\{z_1, z_2, \dots\}$ spans the Hilbert space $\Cal E_P$. 
If, in addition, $z_1, z_2, \dots$ satisfies the linear 
independence condition (1.9),  then it is an orthonormal basis 
for $\Cal E_P$.  
\endroster
These properties are discussed more fully in \cite{4}.  

\proclaim{Definition 1.12}
A {\bf metric operator space} is a pair $(\Cal E, \<\cdot,\cdot\>)$
consisting of a linear supspace $\Cal E\subseteq M=\Cal B(H)$ and an inner 
product $\<\cdot,\cdot\>: \Cal E\times\Cal E\to \Bbb C$ with 
respect to which $\Cal E$ is a separable Hilbert space with the 
following property: for any orthonormal basis $v_1, v_2, \dots$ 
for $\Cal E$ we have 
$$
\|v_1^*\xi\|^2 + \|v_2^*\xi\|^2 + \dots < \infty
$$
for every $\xi\in H$.  
\endproclaim
It is apparent from the preceding remarks 
that the positive operator defined by the 
sum $\sum_kv_kv_k^*$ does not depend on the choice of 
orthonormal basis $(v_k)$, and in fact we can 
associate with $\Cal E$ a unique normal 
completely positive linear map $P_\Cal E$ on $M$ by
$$
P_\Cal E(x) = \sum _k v_kxv_k^*, \qquad x\in M.  
$$
The preceding remarks can now be summarized as follows:
In the von Neumann algebra
$M=\Cal B(H)$, the association 
$\Cal E \leftrightarrow P_\Cal E$ defines a bijective 
correspondence between the set of metric operator spaces 
contained in $M$ and the set of normal 
completely positive linear maps in $\Cal L(M)$.  

\remark{Remark 1.13}
Every normal completely positive linear map can be 
decomposed into a sum of the form 
$P_{\Cal E_0} + c\cdot\iota_M$, where $c$ is a nonnegative 
scalar, $\iota_M4$ is the identity map of $M$, 
and $\Cal E_0$ is a metric operator space 
satisfying the condition 
$\Cal E_0\cap \Bbb C \bold 1 = \{0\}$.  To see that, 
suppose that $\Cal E$ contains the identity operator 
and we set 
$\Cal E_0 = \{v\in \Cal E: \<v, \bold 1\>_\Cal E = 0\}$.  
Then $\Cal E_0$ is a metric operator space satisfying 
$\Cal E_0\cap \Bbb C\bold 1 = \{0\}$, and notice that there
is a positive scalar $c$ such that 
$$
P_\Cal E(x) = P_{\Cal E_0}(x) + cx.  
$$
Indeed, we can choose an orthonormal basis 
$v_0, v_1, v_2, \dots$ for $\Cal E$ so that 
$v_0 = \lambda\bold 1$ is a multiple of $\bold 1$.  
Then $v_1, v_2, \dots$ is an orthonormal basis for 
$\Cal E_0$ and we have 
$$
P_\Cal E(x) = \sum_{k=0}^\infty v_kxv_k^* = 
\sum_{k=1}^\infty v_kxv_k^* + v_0xv_0^* = 
P_{\Cal E_0}(x) + |\lambda|^2x.  
$$
\endremark

\remark{Remark}
While a normal completely positive map $P$ on $M$ 
determines its metric operator space $\Cal E$ uniquely, 
that is not the case for the symbol of $P$.  More 
precisely, if $\Cal E_1$ and $\Cal E_2$ are two metric 
operator spaces with respective completely positive
maps $P_1$ and $P_2$, then $\sigma_{P_1} = \sigma_{P_2}$ 
iff $\Cal E_1 + \Bbb C\bold 1 = \Cal E_2 + \Bbb C\bold 1$.  
We will not make use of that fact, but we do require the 
following more explicit result from which it is  
easily deduced (see Theorem 3.3 and 
Remark 3.18 below for general 
results related to this issue).  
\endremark  

\proclaim{Theorem 1.14}
Let $\Cal E$ be a metric operator space satisfying 
$\Cal E\cap \Bbb C\cdot\bold 1 = \{0\}$, and let 
$P = P_\Cal E$ be its completely positive map.  The 
most general normal completely positive linear map
$Q \in \Cal L(M)$ satisfying $\sigma_Q = \sigma_P$ has
the form
$$
Q(x) = \sum_k(v_k+\lambda_k\bold 1)x(v_k+\lambda_k\bold 1)^* 
+ cx, 
$$
where $v_1, v_2, \dots$ is an orthonormal basis for 
$\Cal E$,  $(\lambda_1, \lambda_2, \dots)$ belongs to $\ell^2$,
and $c$ is a nonnegative scalar.  
\endproclaim
\demo{proof}
In view of (Remark 1.13), it suffices to prove the following 
assertion.  Let $\Cal E$, $\Cal {\tilde E}$ be metric 
operator spaces such that 
$\Cal E\cap \Bbb C\bold 1 = 
\Cal {\tilde E}\cap \Bbb C\bold 1 = \{0\}$, for which
$P_\Cal E$ and $P_{\Cal {\tilde E}}$ have the same 
symbol.  Then there is an orthonormal basis 
$v_1, v_2, \dots$ for $\Cal E$ and an $\ell^2$ sequence
$\lambda_1, \lambda_2, \dots$ such that 
$v_1 + \lambda_1\bold 1, v_2 + \lambda_2\bold 1, \dots$ 
is an orthonormal basis for $\Cal{\tilde E}$.  

To prove this, choose an orthonormal basis 
$u_1, u_2, \dots$ for $\Cal E$ and let $n = \dim\Cal E$.  
If we let $\pi$ be the diagonal representation of 
$M = \Cal B(H)$ on $H^n$
$$
\pi(x) = x\oplus x\oplus \dots, 
$$
then we have 
$$
P_\Cal E(x) = \sum_ku_kxu_k^* = V^*\pi(x)V, 
$$
where $V\in \Cal B(H,H^n)$ is the operator 
$$
V\xi = (u_1^*\xi, u_2^*\xi, \dots).  
$$
Because of (1.9) we have 
$H^n = [\pi(x)\xi: x\in M, \xi \in H]$.  However in this case, 
since $\Cal E\cap \Bbb C\bold 1 = \{0\}$, we claim 
that in fact 
$$
H^n = [(Vx - \pi(x)V)\xi: x\in M, \xi\in H].  \tag{1.15}
$$
To prove that let $K$ be the subspace of $H^n$ 
defined by the right side of (1.15).  Since the 
map $x\in M \mapsto D(x) = Vx - \pi(x)V \in \Cal B(H,K)$
satisfies 
$$
D(xy) = \pi(x)D(y) + D(x)y,
$$
it follows that $K = [D(M)H]$ is invariant under 
operators in $\pi(M)$, and hence the projection 
$p$ onto the orthocomplement of $K$ belongs to 
the commutant of $\pi(M)$ and satisfies 
$pD(x) = 0$ for every $x\in M$.  We have to show 
that $p=0$.  Considering the form
of $\pi$ we find that $p$ is a matrix of scalar operators
$p = (\lambda_{ij}\bold 1)$.  Each row of the scalar 
matrix $(\lambda_{ij})$ belongs to $\ell^2$, and the 
condition $pD(x) = 0$ implies that for every 
$i = 1,2,\dots$ and every $x\in M$ we have 
$$
\sum_j\lambda_{ij}(u_j^*x - xu_j^*) = 0.  
$$
Thus for every $i$ the operator 
$w_i = \sum_j\bar\lambda_{ij}u_j$ is an element of 
$\Cal E$ which commutes with every operator in M, and 
is therefore a scalar multiple of the identity.  
Since $\Cal E\cap \Bbb C\bold 1 = \{0\}$ we conclude
that $w_1 = w_2 = \dots = 0$.  Hence by (1.9) the 
matrix $p = (\lambda_{ij}\bold 1)$ is zero, proving 
(1.15).  

Now let $m = \dim \Cal{\tilde E}$ and let 
$$
\tilde\pi(x) = x\oplus x\oplus \dots 
$$
be the corresponding representation of $M$ on 
$H^m$.  Let $\tilde u_1, \tilde u_2, \dots$ be an 
orthonormal basis for $\Cal{\tilde E}$ and let 
$\tilde V: H\to H^m$ be the associated operator 
$$
\tilde V\xi = (\tilde u_1^*\xi, \tilde u_2^*\xi, \dots).  
$$
Then there is a corresponding derivation 
$\tilde D: M\to \Cal B(H,K)$ which is defined by
$\tilde D(x) = \tilde Vx - \tilde\pi(x)\tilde V$.  
We claim next that there is a (necessarily unique) 
unitary operator $W: H^n \to H^m$ satisfying 
$$
WD(x) = \tilde D(x), \qquad x\in M.  \tag{1.16}
$$
Noting that $H^n = [D(M)H]$ and $H^m = [\tilde D(M)H]$, 
it is clearly enough to show that for all $\xi, \eta\in H$
we have 
$$
\<D(x)\xi, D(y)\eta\> = \<\tilde D(x)\xi, \tilde D(y)\eta\>.  
\tag{1.17}
$$
For that, we write 
$$
\align
D(y)^*D(x) &= (Vy - \pi(y)V)^*(Vx - \pi(x)V) \\
&= V^*\pi(xy)V - y^*V^*\pi(x)V - V^*\pi(y^*)Vx + y^*V^*Vx\\
&= \sigma_{P_\Cal E}(dy^*\, dx) .  
\endalign
$$
By hypothesis the symbols of $P_\Cal E$ and $P_{\tilde\Cal E}$
agree, hence the right side is 
$$
\sigma_{\tilde{\Cal E}}(dy^*\,dx) = 
\tilde D(y)^*\tilde D(x), 
$$
and formula (1.17) follows.  

Note that $W\pi(x) = \tilde \pi(x)W$ for every $x$.  
Indeed, fixing $x$ and choosing a vector in $H^n$ 
of the form $\eta = D(y)\xi$ for $y\in M$, $\xi\in H$ 
we have 
$$
\align
W\pi(x)\eta &= W\pi(x)D(y)\xi = WD(xy)\xi - WD(x)y\xi \\
&= \tilde D(xy)\xi - \tilde D(x)y\xi = \tilde\pi(x)\tilde D(y)\xi 
= \tilde\pi(x)WD(y)\xi.  
\endalign
$$
The assertion follows because $H^n$ is spanned by such 
vectors $\eta$.  

In particular, $\pi$ and $\tilde\pi$ are equivalent 
representations of $M$.  Hence $m=n$ and therefore $H^m=H^n$.  
Moreover, $W$ belongs to the commutant of $\pi$ and 
hence there is a unitary matrix $(\lambda_{ij})$ of 
complex scalars such that $W = (\lambda_{ij}\bold 1)$.  
If we now look at the components of the operator equation 
$WD(x) = \tilde D(x)$ we find that for every $i = 1, 2, \dots$ 
$$
\sum_j\lambda_{ij}(u_j^*x - xu_j^*) 
= \tilde u_i^*x - x\tilde u_i^*.  
$$
Thus we can define a new orthonormal basis $v_1, v_2, \dots $
for $\Cal E$ by 
$$
v_i = \sum_j \bar\lambda_{ij}u_j.  
$$
The preceding equation relating the $u_j$ to the $\tilde u_j$ 
becomes
$$
v_i^*x - xv_i^* = \tilde u_i^*x - x \tilde u_i^*
$$
for every operator $x\in M=\Cal B(H)$.  It follows that 
the operators $\tilde u_i - v_i$ commute with all bounded 
operators and therefore must be scalar multiples 
of the identity operator, say
$$
\tilde u_i = v_i + \lambda_i\bold 1.  \tag{1.18}
$$
The fact that both series $\sum_j v_jv_j^*$ and 
$\sum \tilde u_j\tilde u_j^*$ converge strongly implies 
that $\sum_j|\lambda|^2 < \infty$, and thus (1.18) provides
an orthonormal basis $\tilde u_1, \tilde u_2, \dots$ 
for $\Cal{\tilde E}$ of the required form.
\qed
\enddemo

\proclaim{Corollary}
Let $\Cal E$ and $P = P_\Cal E$ satisfy the hypotheses 
of Theorem 1.14, and let $k$ be an operator in $M$.  
The following are equivalent:
\roster
\item
$Q(x) = P(x) +kx + xk^*$ is completely positive.  
\item
$k$ has a (necessarily unique) decompostion of the form
$k = v + c\bold 1$, where $v\in \Cal E$ and $c$ is a 
complex number satisfying 
$c+\bar c\geq \<v,v\>_\Cal E$.
\endroster
\endproclaim
\demo{proof of (2)$\implies$(1)}
Let $k=v + 1/2\<v,v\>_\Cal E + d\bold 1$ where $v\in \Cal E$
and $d$ is a complex number with nonnegative real part.  
Choose an orthonormal basis $v_1, v_2,\dots$ for $\Cal E$ 
and set $\lambda_k = \<v,v_k\>$.  Then for 
$\alpha = d + \bar d - \<v,v\>_\Cal E \geq 0$ we have 
$$
kx + xk^* = vx + xv^* + (\<v,v\>_\Cal E + \alpha) x = 
\sum_k(\lambda_kv_kx + \bar\lambda_k xv_k^* + 
|\lambda_k|^2 x) + \alpha x,
$$
hence the map $Q$ of (1) can be written 
$$
Q(x) = \sum_k(v_k + \bar\lambda_k\bold 1)x(v_k+\bar\lambda_k\bold 1)^*
+ \alpha x
$$
which is obviously completely positive.  
\enddemo

Before proving the opposite implication we collect 
two elementary observations.  

\proclaim{Lemma 1.19}
\roster
\item
The only completely positive linear map of $M = \Cal B(H)$ 
having symbol $0$ is of the form $L(x) = cx$ where 
$c$ is a nonnegative scalar.  
\item
The only operator $k\in M$ for which $P(x) = kx + xk^*$ is 
a completely positive map is of the form $k = z\bold 1$, 
where $z$ is a complex number having nonnegative real part.  
\endroster
\endproclaim
\demo{proof}
Suppose that $P$ is a completely positive linear map on 
$M$ for which $\sigma_P = 0$.  Let $P(x) = V^*\pi(x)V$ be 
a Stinespring representation for $P$.  From the definition
of $\sigma_P$ one finds that 
$$
\sigma_P(dx\,dy) = (V^*\pi(x) - xV^*)(\pi(y)V - Vy).  
$$
Therefore 
$(\pi(x)V - Vx)^*(\pi(x)V - Vx) = \sigma_P(d(x^*)\,dx) = 0$, 
hence $Vx = \pi(x)V$ for every $x$.  It follows that 
$V^*Vx = xV^*V$ for all $x\in M = \Cal B(H)$ so that there is
a nonnegative scalar $x$ such that $V^*V = c\bold 1$.  
Thus $P(x) = V^*\pi(x)V = V^*Vx = cx$ as asserted.  

For the second assertion, suppose that $x\mapsto kx + xk^*$ 
is completely positive.  The symbol of this map vanishes, 
so by what was just proved there is a nonnegative scalar 
$c$ such that $kx + xk^* = cx$ for every $x$.  Taking $x$ 
to be an arbitrary projection $p$ we find that 
$(\bold 1-p)kp = 0$, hence $k$ must be a scalar 
$k = \lambda \bold 1$.  The formula $kx + xk^* = cx$ 
implies that $\lambda + \bar\lambda = c \geq 0$, 
and the Lemma is proved.   
\qed
\enddemo

\demo{proof of (1)$\implies$(2)}
Suppose that $k$ is an operator in $M$ such that 
$Q(x) = P(x) +kx + xk^*$ is completely positive.  Then 
$Q$ and $P$ have the same symbol, hence by Theorem 1.14 there
is an $\ell^2$ sequence $\lambda = (\lambda_1, \lambda_2, \dots)$
and an orthonormal basis $v_1, v_2, \dots$ for $\Cal E$ such 
that 
$$
Q(x) = \sum_k(v_k + \lambda_k\bold 1)x(v_k + \lambda_k\bold 1)^*.  
$$
Define an element $v\in \Cal E$ by $v = \sum_k\bar\lambda_kv_k$, 
and write $|v|^2 = \<v,v\>_\Cal E = \sum_k|\lambda_k|^2$.  We find 
that 
$$
\align
Q(x) &= \sum_kv_kxv_k^* + vx + xv^* + |v|^2 x \\
&= P(x) + (v + 1/2|v|^2\bold 1)x + x(v + 1/2|v|^2\bold 1)^* \\
&= P(x) + kx + xk^*.  
\endalign
$$
It follows that the operator $\ell = k - v - 1/2|v|^2\bold 1$ 
has the property that $x \mapsto \ell x + x \ell^*$ is a 
completely positive linear map on $M = \Cal B(H)$.  From part
(2) of the preceding Lemma we find that there is a complex scalar
$z$ having nonnegative real part such that 
$\ell = z\bold 1$, and the required representation
$$
k = v + (1/2|v|^2 + z)\bold 1
$$
follows\qed
\enddemo

We conclude this section by reformulating
a known description of 
the bounded generators of semigroups of completely 
positive maps in terms of metric operator spaces.

\proclaim{Proposition 1.20}
Let $L\in \Cal L(M)$ be an operator which generates 
a semigroup $P_t = \exp(tL)$, $t\geq 0$ of normal completely
positive maps on $M$.  Then there is a metric operator space 
$\Cal E$ satisfying $\Cal E\cap \Bbb C\bold 1 = \{0\}$, 
and an operator $z\in M$ such that 
$$
L(x) = P_\Cal E(x) + zx + xz^* \qquad x\in M.  
$$   
\endproclaim
\demo{proof}
Using a general result of Christensen and Evans \cite{7}, one
can find a completely positive linear map $Q: M\to M$ and 
an element $z\in M$ such that
$$
L(x) = Q(x) + zx + xz^*, \qquad x\in M.  \tag{1.21}
$$
Since $L$ is a bounded operator that generates a semigroup
of normal completely positive maps, it must itself be a normal
linear map on $M$; and we may conclude from 
(1.21) that $Q$ is normal.  Let $\Cal E$ be a metric operator
space such that $Q = P_\Cal E$.  The remark following Definition 
1.12 shows that we can arrange $\Cal E\cap\Bbb C\bold 1=\{0\}$
by adjusting $z$ if necessary
\qed
\enddemo

\remark{Remarks}
Unlike the case of completely positive maps, the correspondence
between metric operator spaces and generators of CP semigroups 
is not quite one-to-one.  However, if $(\Cal E_1,z_1)$ and 
$(\Cal E_2,z_2)$ are two pairs which serve to represent a given
generator $L$ as in Proposition 1.20 and which satisfy 
$$
\Cal E_1\cap \Bbb C\bold 1 = 
\Cal E_2\cap\Bbb C\bold 1 = \{0\}, \tag{1.22}
$$
then it follows from Theorem 1.14 that 
$\Cal E_1+\Bbb C\bold 1 = \Cal E_2+\Bbb C\bold 1$, and hence
$\dim \Cal E_1 = \dim \Cal E_2$.  Thus we can make the following 

\proclaim{Definition 1.23}
Let $L$ be a bounded operator on $\Cal B(H)$ which generates
a semigroup of normal completely positive maps on 
$\Cal B(H)$.  The {\bf rank} of $L$ is defined as the dimension of 
any metric operator space $\Cal E$ for which 
$\Cal E\cap\Bbb C\bold 1 = \{0\}$ and which gives rise
to a representation of $L$ in the form 
$$
L(x) = P_\Cal E(x) + kx + xk^*, \qquad x\in \Cal B(H)
$$
where $k$ is some operator in $\Cal B(H)$.  
\endproclaim

We also point out that, with a little more care, one 
can recover the inner product on $\Cal E$ from the 
properties of $L$ (see Remark 3.18 below).

\subheading{2. Units and covariance function}

In \cite{4}, a notion of index 
for CP semigroups $P$ was introduced 
which directly generalizes the definition of numerical index
of \esg s.  Briefly, a {\bf unit} of $P$ is a strongly 
continuous semigroup $T = \{T(t): t\geq 0\}$ of bounded 
operators in $M$ for which there is a real constant 
$k$ with the property that for every $t\geq 0$, the operator
mapping 
$$
x\in M \mapsto e^{kt}P_t(x) - T(t)xT(t)^* 
$$
is completely positive.  Let $\Cal U_P$ be the set of all 
units of $P$.  It is possible for $\Cal U_P$ to be the 
empty set; indeed there are \esg s with this 
property \cite{9,10}.  
But if $\Cal U_P\neq \emptyset$ then one 
can define a function 
$$
c_P: \Cal U_P\times \Cal U_P \to \Bbb C, 
$$
called the {\bf covariance function} of the semigroup $P$, 
as follows.  Note first that for every positive 
$t$, there is a unique metric operator space associated
with the completely positive map $P_t$; we will write
$\Cal E_P(t)$ for this metric operator space.  The 
most elementary properties of the family 
$$
\Cal E_P = \{\Cal E_P(t): t > 0\}
$$
are as follows:
\roster
\item"{(2.1.1)}" 
Each $\Cal E_P(t)$ is a separable Hilbert space.
\item"{(2.1.2)}"
$\Cal E_P(s+t)$ is spanned as a Hilbert space
by the set of products 
$$
\{xy: x\in \Cal E_P(s), y\in \Cal E_P(t)\}.
$$  
\endroster

\remark{Remarks}
Regarding (2.1.2), it is shown in \cite{4, Theorem 1.12} that 
operator multiplication 
$$
u\otimes v\in\Cal E_P(s)\otimes\Cal E_P(t)
\mapsto uv\in\Cal E_P(s+t)
$$ 
extends uniquely to a bounded linear 
operator from the Hilbert space 
$\Cal E_P(s)\otimes\Cal E_P(t)$ to $\Cal E_P(s+t)$ whose
adjoint is an isometry from $\Cal E_P(s+t)$ onto a closed
subspace of $\Cal E_P(s)\otimes \Cal E_P(t)$.  
\endremark

In order to define the covariance function, choose 
$T_1, T_2\in \Cal U_P$, and 
fix $t>0$.  The condition (2.1) implies that both operators
$T_1(t)$ and $T_2(t)$ belong to $\Cal E_P(t)$ and 
thus we can form their inner product 
$\<T_1(t), T_2(t)\>_{\Cal E_P(t)}$ as elements of this Hilbert space.  
More generally, if we are given an arbitrary finite 
partition 
$$
\Cal P = \{0 = t_0 < t_1 < \dots < t_n = t\}
$$
of the interval $[0,t]$ then we define a function 
$f_{\Cal P,t}: \Cal U_p\times\Cal U_P \to\Bbb C$ as follows
$$
f_{\Cal P,t}(T_1,T_2) = 
\prod_{k=1}^n 
\<T_1((t_k-t_{k-1}),T_2(t_k-t_{k-1})\>_{\Cal E_P(t_k-t_{k-1})}.  
$$
It was shown in \cite{4} that there 
is a unique complex number $c_P(T_1,T_2)$ such that 
for every $t>0$
$$
\lim_\Cal P f_{\Cal P,t}(T_1,T_2) = e^{tc_P(T_1,T_2)}, \tag{2.2}
$$
the limit being taken over the increasing directed set 
of all finite partitions $\Cal P$ of $[0,t]$.  
That defines the covariance
function $c_P: \Cal U_P\times \Cal U_P \to \Bbb C$ of any
CP semigroup $P$ for which $\Cal U_P\neq \emptyset$.  

The covariance function is conditionally positive definite 
in the sense that if $T_1, T_2,\dots,T_n\in\Cal U_P$ and 
$\lambda_1, \lambda_2,\dots,\lambda_n\in \Bbb C$ satisfy 
$\lambda_1+\lambda_2+\dots+\lambda_n=0$, then 
$$
\sum_{i,j=1}^n\lambda_i\bar\lambda_jc_P(T_i,T_j)\geq 0
$$
(see \cite{cpindex, Proposition 2.7}).  More generally, if 
we are given any nonempty set $X$ and a conditionally positive 
definite function $c: X\times X\to\Bbb C$ then there is a 
natural way to construct a Hilbert space $H(X,c)$.  Briefly, 
$c$ defines a positive semidefinite sesquilinear form 
$\<\cdot,\cdot\>$ on the vector space $V$ of all finitely 
nonzero functions $f: X\to \Bbb C$ for which 
$$
\sum_{x\in X}f(x)=0,
$$
by way of 
$$
\<f,g\>=\sum_{x,y\in X}f(x)\bar g(y)c(x,y), 
$$
and $H(X,c)$ is obtained by completing the inner product 
space obtained by promoting $\<\cdot,\cdot\>$ to the 
quotient of $V$ by the subspace 
$$
N=\{f\in V: \<f,f\>=0\}.  
$$

The index $d_*(P)$ of a CP semigroup $P$ is defined by
$$
d_*(P)=\dim H(\Cal U_P,c_P)
$$
in the case where $\Cal U_P\neq \emptyset$, and is 
defined by $d_*(P)=2^{\aleph_0}$ if $\Cal U_P=\emptyset$.  
Notice that in order to calculate $d_*(P)$ one must 
calculate a) the set $\Cal U_P$ of all units of $P$, and 
b) the covariance function $c_P: \Cal U_P\times\Cal U_P\to\Bbb C$; 
moreover, this must be done explicitly enough so that the dimension
of $H(\Cal U_P,c_P)$ is apparent.  The purpose of this section 
is to carry out these calculations for the case of CP 
semigroups with {\it bounded} generator in terms of the 
structures associated with the generator by Proposition 1.20.  
The principal result is Theorem 2.3 below.

\remark{Remarks}
The covariance function $c_P:\Cal U_P\times\Cal U_P\to\Bbb C$ is
a conditionally positive definite function having the 
property that for every $t>0$, 
$$
T_1, T_2\in \Cal U_P \mapsto e^{tc_P(T_1,T_2)} - 
\<T_1(t),T_2(t)\>_{\Cal E_P(t)}
$$
is a positive definite function.  It is a simple exercise 
to show that if $d: \Cal U_P\times\Cal U_P\to\Bbb C$ is any 
other function for which 
$$
T_1, T_2\in \Cal U_P \mapsto e^{td(T_1,T_2)} - 
\<T_1(t),T_2(t)\>_{\Cal E_P(t)}
$$
is positive definite for every $t>0$, then the difference 
$d-c_P$ is a positive definite function on $\Cal U_P\times\Cal U_P$.  
Thus the covariance function $c_P$ is characterized in this 
sense as the ``smallest" function $c:\Cal U_P\times\Cal U_P\to\Bbb C$ 
with the property that 
$e^{tc(T_1,T_2)}$ dominates the inner products 
$\<T_1(t),T_2(t)\>_{\Cal E_P(t)}$ for every $t>0$.  
\endremark

\proclaim{Theorem 2.3}
Let $L\in \Cal L(M)$ be a bounded operator which generates
a CP semigroup on $M = \Cal B(H)$.  
Let $\Cal E$ be a metric
operator space satisfying the conditions of Proposition 1.20, 
so that $L$ has the form
$$
L(x) = P_\Cal E(x) + kx + xk^* \tag{2.3.1}
$$
for some $k\in M$.  The units of the CP semigroup 
$P = \{\exp(tL): t\geq 0\}$ are described 
in terms of $\Cal E$ and $k$ as follows.  
For every $(c,v)\in \Bbb C\times \Cal E$, let 
$T_{(c,v)}$ be the operator semigroup 
$$
T_{(c,v)}(t) = e^{ct}\exp t(v + k), \qquad t\geq 0.  \tag{2.3.2}
$$
Then $T_{(c,v)}$ is a unit of $P$ and 
the map $(c,v)\in \Bbb C\times \Cal E \mapsto T_{(c,v)}$ is 
a bijection of $\Bbb C\times \Cal E$ onto the set 
$\Cal U_P$ of units of $P$.  

The covariance function 
$c_P: \Cal U_P\times \Cal U_P\mapsto \Bbb C$ of $P$ is given by
$$
c_P(T_{(c_1,v_1)},T_{(c_2,v_2)}) = c_1 + \bar c_2 + \<v_1,v_2\>_\Cal E,  
\tag{2.3.2}
$$
and the index of $P$ is 
$d_*(P) = \dim \Cal E$.
\endproclaim

\demo{proof}
Notice first that the map $(c,v)\to T_{(c,v)}$ 
from $\Bbb C\times \Cal E$ to operator semigroups is one-to-one.  
Indeed, choosing complex numbers $c_1, c_2$ and elements 
$v_1, v_2\in \Cal E$ such that $T_{(c_1,v_1)}(t) = T_{(c_2,v_2)}(t)$
for every $t$, it follows that the generators of these two semigroups
are equal.  Using (2.3.2) we find that
$$
c_1\bold 1 + v_1 + k = c_2\bold 1 + v_2 + k.  
$$
Cancelling $k$ and using  $\Cal E\cap \Bbb C\bold 1 = \{0\}$, 
we obtain $v_1 = v_2$ and $c_1 = c_2$ as asserted.  

Now fix $(c,v)\in \Bbb C\times \Cal E$.  We will show that
$T_{(c,v)}$ is a unit of $P$.  In order to do that, we must 
find a real constant $\alpha$ such that each mapping
$$
x\in M \mapsto e^{t\alpha}P_t(x) - T_{(c,v)}(t)xT_{(c,v)}(t)^* \tag{2.4}
$$ 
is completely positive, $t\geq 0$.  Noting that 
$T_{(c,v)}(t) = e^{ct}T_{(0,v)}(t)$, it is clearly enough to prove 
(2.4) for $c=0$; in that case we will show that 
(2.4) is true for the constant $\alpha = \<v,v\>_\Cal E$.  

To that end, we set 
$$
L_0(x) = (v+k)x + x(v+k)^*,
$$
and we claim first that for $\alpha = \<v,v\>_\Cal E$ we have 
$$
L_0 \leq L + \alpha\cdot \iota_M, \tag{2.5}
$$
$\iota_M$ denoting the identity map of $M$.  Equivalently, 
after choosing an orthonormal basis $v_1, v_2, \dots $ for 
$\Cal E$, we want to show that the mapping 
$$
x \mapsto \sum_j v_jxv_j^* - vx -xv^* + \alpha x  \tag{2.6}
$$
is completely positive.  Indeed, since $v\in \Cal E$ and 
$(v_j)$ is an orthonormal basis, we can
find a sequence $\lambda\in \ell^2$ such that 
$v = \sum_j\bar\lambda_jv_j$ (note the complex conjugate).  
Then
$$
\alpha = \<v,v\> = \sum_j|\lambda_j|^2, 
$$
and hence the term on the right side of (2.6) 
can be collected as follows, 
$$
\sum_j(v_j - \lambda_j\bold 1)x(v_j - \lambda_j)^*.  
$$
The latter formula obviously defines a completely positive map.  

In order to pass from (2.5) to its exponentiated 
version (2.4), we require
\proclaim{Lemma 2.7}
Let $L_1$, $L_2$ belong to $\Cal L(M)$.  Suppose that 
both generate CP semigroups and that $L_2 - L_1$ is completely
positive.  Then for every $t\geq 0$ the map 
$$
\exp(tL_2) - \exp(tL_1)
$$
is completely positive.  
\endproclaim
\demo{proof}
Since the hypotheses on $L_1$ and $L_2$ are invariant under 
scaling by positive constants, it is enough to prove the 
assertion for $t=1$.  We can write
$L_2 = L_1 + R$ where $R$ is a completely positive map.  
By the Lie product formula \cite{12, page 245}, 
we have 
$$
\exp(L_2) = \exp(L_1 + R) = 
\lim_{n\to \infty}(\exp(1/nL_1)\exp(1/nR))^n,
$$
the limit on $n$ existing relative to the operator norm
on $\Cal L(M)$.  Thus it suffices to show that 
for every $n$, we have 
$$
(\exp(1/nL_1)\exp(1/nR))^n \geq \exp L_1.  \tag{2.8}
$$
To see the latter, note that for {\it completely positive} 
maps $A_k, B_k$, $k=1,2$ we have 
$$
B_1 \geq A_1 \text{ and } B_2\geq A_2 
\implies B_1B_2\geq A_1A_2.  
$$
Indeed, this follows from the fact that a composition 
of completely positive maps is completely positive, so 
that $A_1(B_2-A_2) \geq 0$ and $(B_1 - A_1)B_2 \geq 0$ 
together imply that $B_1B_2 \geq A_1B_2 \geq A_1A_2$, 
and the assertion follows.  

We apply this to (2.8) as follows.  Letting $\iota_M$ denote
the identity map of $M$ we have 
$$
\exp(1/nR) = \iota_M + 1/nR + 1/2(1/nR)^2 + \dots 
\geq \iota_M
$$
because $R$ is completely positive.  Hence 
$$
\exp(1/nL_1) \exp(1/nR) \geq \exp(1/nL_1).  
$$
For the same reason, 
$$
(\exp(1/nL_1) \exp(1/nR))^2 \geq (\exp(1/nL_1))^2,   
$$
and so on until we obtain 
$$
(\exp(1/nL_1) \exp(1/nR))^n \geq (\exp(1/nL_1))^n = \exp(L_1).   
$$
This establishes (2.8) and completes the proof of Lemma 2.7\qed
\enddemo

Applying Lemma 2.7 to (2.5) we obtain
$$
\exp(tL_0) \leq \exp(t(L + \alpha\iota_M)) = e^{t\alpha}P_t.  
$$
Noting that $L_0(x) = (v + k)x + x(v + k)^*$ is the generator
of the semigroup 
$$
x\mapsto T_{(0,v)}(t)xT_{(0,v)}(t)^*,
$$ 
we obtain formula (2.4) for the case $c=0$.  This completes 
the proof that each operator semigroup 
of the form $T_{(c,v)}\in \Cal U_P$ 
is a unit of $P$.  

We show next that the map
$(c,v)\in \Bbb R\times \Cal E \mapsto T_{(c,v)}\in\Cal U_P$.  
is surjective.  For that,
the following result is essential.  

\proclaim{Lemma 2.9}Let $L$ be a bounded operator 
in $\Cal L(M)$ which generates a unital CP semigroup.  
Then every semigroup $T\in \Cal U_P$ has a bounded 
generator.  
\endproclaim
\demo{proof}Let $T$ be an operator semigroup with the property
that for every $t\geq 0$ the mapping 
$$
x\mapsto e^{\alpha t}P_t(x) - T(t)xT(t)^* \tag{2.10}
$$
is completely positive, 
$P_t$ denotine $\exp(tL)$ and $\alpha$ being some 
real constant.  To show that the generator of $T$ is a bounded 
operator, it is enough to show that $T$ is continuous 
relative to the operator norm on $M$
$$
\lim_{t\to 0+}\|T(t) - \bold 1\| = 0.  \tag{2.11}
$$
To prove (2.11) we make use of the symbol as follows.  
Fix $t > 0$ and consider the operator mapping 
$$
L_t(x) = T(t)xT(t)^*, \qquad x\in M.  
$$
The symbol of $L_t$ is found to be
$$
\sigma_{L_t}(dx\, dy) = (T(t)x - xT(t))(T(t)^*y - yT(t)^*).  
$$
Now the symbol of a completely positive map $R$ satisfies
$\sigma_R(dx\,dx^*) \geq 0$.  Hence 
if $R_1$ and $R_2$ are completely positive maps satisfying
$0 \leq R_1 \leq R_2$ then we have 
$0 \leq \sigma_{L_1}(dx\,dx^*) \leq \sigma_{L_2}(dx\,dx^*)$.  
Thus, (2,10) implies that for all $x\in M$,
$$
\align
(T(t)x - xT(t))(T(t)x - xT(t))^* &= \sigma_{L_t}(dx\,dx^*) 
\leq e^{t\alpha}\sigma_{P_t}(dx\,dx^*) \\
&= e^{t\alpha}\sigma_{P_t-\iota_M}(dx\,dx^*), 
\endalign
$$
the last equality resulting from the fact that the identity
map $\iota_M$ of $M$ has symbol zero.  
From (1.4) and the previous formula 
we conclude that for all $x\in M$ satisfying 
$\|x\| \leq 1$, we have 
$$
\|T(t)x - xT(t)\|^2 
= e^{t\alpha}\|\sigma_{P_t - \iota_M}(dx\,dx^*)\| \leq 
4 e^{t\alpha}\|P_t-\iota_M\|, 
$$
$\|P_t-\iota_M\|$ denoting the norm of $P_t-\iota_M$ as 
an element of $\Cal L(M)$.  Now 
$$
\|P_t-\iota_M\| = \|\exp(tL) - \iota_M\| \to 0
$$
as $t\to 0$ because $L$ is bounded.  Since the norm of 
a derivation of $M = \Cal B(H)$ of the form 
$D(x) = Tx - xT$ satisfies inequalities of the form 
$$
\inf_{\lambda\in \Bbb C}\|T - \lambda\bold 1\| \leq \|D\| 
\leq 2 \inf_{\lambda\in \Bbb C}\|T - \lambda\bold 1\|,  
$$ 
it follows that 
$$
\inf_{\lambda\in \Bbb C}\|T(t) - \lambda\bold 1\| \leq 
\sup_{\|x\|\leq 1}\|T(t)x - xT(t)\| \to 0, 
$$
as $t\to 0$.  Thus there exist complex scalars $\lambda_t$ 
such that $\|T(t) - \lambda_t\bold 1\|\to 0$ as $t\to 0$.  
Since the semigroup $\{T(t) : t\geq 0\}$ is strongly continuous,
$T(t)$ must tend to $\bold 1$ in the strong operator topology
as $t\to 0$; hence $\lambda_t \to 1$ as $t\to 0$ and
(2.11) follows \qed
\enddemo

Now choose any unit $T\in \Cal U_P$.  
There is a real constant
$\alpha$ such that for every $t$ the mapping 
$$
x\mapsto e^{t\alpha}P_t(x) - T(t)xT(t)^*
$$
is completely positive.  We will show
that there is an element $(c,v)\in \Bbb C\times \Cal E$ 
such that $T = T_{(c,v)}$.  
By replacing $T$ with the semigroup
$\{e^{-\alpha t/2}T(t): t\geq 0\}$ (and adjusting 
$c$ accordingly), 
we may clearly assume that $\alpha = 0$.  
By Lemma 2.9, there 
is a {\it bounded} operator $a\in M$ such that 
$$
T(t) = \exp(ta), \qquad t\geq 0  
$$
and we have to show that $a$ has the form
$$
a = c\bold 1 + v + k \tag{2.12}
$$
for some complex scalar $c$ and some $v\in \Cal E$.  
For that, we claim first that the operator mapping 
$$
R(x) = L(x) -ax - xa^*  \tag{2.13}
$$
is completely positive.  Indeed, since for every $t>0$ the map
$$
x\mapsto P_t(x) -e^{ta}xe^{ta^*} =
 (P_t(x) - x) -(e^{ta}xe^{ta^*} - x)
$$
is completely positive, we may divide the latter by
$t$ and take the limit as $t\to 0+$ to obtain (2.13),
after noting that 
$$
\align
\lim_{t\to 0+}t^{-1}(P_t(x) - x) &= L(x), \qquad \text{and} \\
\lim_{t\to 0+}t^{-1}(e^{ta}xe^{ta^*} - x) &= ax + xa^*.  
\endalign
$$  
Using (2.3.1) we can write $R$ in the form 
$$
R(x) = P_\Cal E(x) + (k - a)x + x(k - a)^*.  
$$
The Corollary of Theorem 1.14 implies that
there is a complex number $d$ and an element $u\in \Cal E$ such 
that 
$$
k - a = d\bold 1 + u, 
$$
and the required representation (2.12) follows after
taking $c = -d$ and $v = -u$.  

It remains to compute the covariance function $c_P$ 
of formula (2.2) in these coordinates.

\proclaim{Lemma 2.14}
Let $v_1, \dots, v_n\in \Cal E$ and let $T_j$ be the 
unit of $P$ defined by 
$$
T_j(t) = \exp{t(v_j + k)}, \qquad t\geq 0, 1\leq j\leq n.  
$$
Then for every $\lambda_1, \dots, \lambda_n\in \Bbb C$ we 
have 
$$
\sum_{i,j}\lambda_i\bar\lambda_j\<v_i, v_j\>_\Cal E \leq 
\sum_{i,j}\lambda_i\bar\lambda_j c_P(T_i, T_j).  
$$
\endproclaim
\demo{proof}
Notice that if an $n$-tuple $\lambda_1,\dots,\lambda_n$ satisfies
the required inequality and $c$ is an arbitrary complex number, 
then so does $c\lambda_1, \dots , c\lambda_n$.  Thus it is 
enough to prove Lemma 2.14 for $n$-tuples $\lambda_k$ satisfying
$\lambda_1 + \lambda_2 + \dots + \lambda_n = 1$.  

Choose  $T_1, \dots, T_n\in\Cal U_P$.  By the remarks following 
(2.2) we have 
$$
\sum_{i,j}\lambda_i\bar\lambda_j\<T_i(t),T_j(t)\>_{\Cal E_P(t)} \leq 
\sum_{i,j}\lambda_i\bar\lambda_je^{tc_P(T_i,T_j)}.  
$$
Equivalently, if for every $t>0$ we set 
$$
A(t) = \lambda_1T_1(t) + \dots +\lambda_nT_n(t), 
$$
then $A(t) \in \Cal E_P(t)$ and we have 
$$
\<A(t), A(t)\>_{\Cal E_P(t)}\leq 
\sum_{i,j}\lambda_i\bar\lambda_j e^{tc_P(T_i,T_j)}.   
$$
It follows that the mapping 
$$
x\mapsto (\sum_{i,j}\lambda_i\bar\lambda_j e^{tc_P(T_i,T_j)})P_t(x)
- A(t)xA(t)^*
$$
is completely positive.  Since $\sum_j\lambda_j = 1$, this implies
that 
$$
x\mapsto \sum_{i,j}\lambda_i\bar\lambda_j( e^{tc_P(T_i,T_j)}P_t(x)-x)
- (A(t)xA(t)^* -x)
$$
is completely positive.  
Notice that $A(t)$ is differentiable at $t=0+$ and that 
$A^\prime(0+) = \sum_j\lambda_jv_j + k$.  
Thus if we divide by $t$ and take $\lim_{t\to 0+}$ 
we find that 
$$
x\mapsto (L(x) + \sum_{i,j}\lambda_i\bar\lambda_jc_P(T_i,T_j)) -
((\sum_j\lambda_jv_j)x + x(\sum_j\lambda_jv_j)^* + kx + xk^*)
$$
is a completely positive map.  Noting that 
$L(x) = P_\Cal E(x) + kx + xk^*$, the terms involving 
$k$ and $k^*$ cancel and we are left with a completely positive 
map of the form
$$
x\mapsto P_\Cal E(x) - v(\lambda)x - 
xv(\lambda)^* + (\sum_{i,j}\lambda_i\bar\lambda_jc_P(T_i,T_j))x.  
$$
where $v(\lambda) = \sum_j\lambda_jv_j\in \Cal E$.  
From the corollary of Theorem 1.14 we deduce the required inequality
$$
\sum_{i,j}\lambda_i\bar\lambda_jc_P(T_i,T_j) \geq 
\<v(\lambda), v(\lambda)\>_\Cal E = 
\sum_{i,j}\lambda_i\bar\lambda_j\<v_i,v_j\>_\Cal E.
$$
\qed
\enddemo

\proclaim{Lemma 2.15}
Let $v\in \Cal E$ and let $T\in \Cal U_P$ be the semigroup 
$T(t) = \exp{t(v+k)}$, $t\geq 0$. Then 
$$
c_P(T,T) \leq \<v,v\>_\Cal E.  
$$
\endproclaim
\demo{proof}
It has already been shown (see (2.4)) that 
the mapping 
$$
x\mapsto e^{t\<v,v\>_\Cal E}P_t(x) - T(t)xT(t)^*
$$
is completely positive.  It follows from the definition
of $\Cal E_P(t)$ that 
$$
\<T(t),T(t)\>_{\Cal E_P(t)} \leq e^{t\<v,v\>_\Cal E}
$$
for every $t>0$.  Thus for every finite partition 
$$
\Cal P = \{0 = t_0 < t_1 < \dots < t_n=t\}
$$
of the interval $[0,t]$ we have 
$$
\prod_{k=1}^n\<T(t_k - t_{k-1}),T(t_k-t_{k-1})\>
_{\Cal E_P(t_k-t_{k-1})} \leq 
e^{t\<v,v\>_\Cal E}.  
$$
Noting the definition (2.2) of $c_P$ we conclude
that 
$$
e^{tc_P(T,T)} \leq e^{t\<v,v\>_\Cal E}
$$
for every $t>0$, and the asserted inequality follows
\qed
\enddemo

To complete the proof of Theorem 2.3, choose complex 
numbers $c_1, c_2$, choose  $v_1, v_2\in \Cal E$, and let 
$T_1$,
$T_2$ be the units of $P$ defined by 
$$
T_k(t) = T_{(0,v_k)}e^{ct}\exp{t(v_k + k)}, \qquad t\geq 0.  
$$
Consider the self-adjoint $2\times 2$ matrices $A = (a_{ij})$ 
and $B = (b_{ij})$ defined by 
$$
\align
a_{ij} &= \<v_i,v_j\>_\Cal E, \\
b_{ij} &= c_P(T_i,T_j).  
\endalign
$$
Lemma 2.14 implies that $B - A\geq 0$, while Lemma 2.15 
implies that the diagonal terms of $B - A$ are nonpositive.  
Hence the trace of $B - A$ is nonpositive and it follows 
that $A = B$.  Comparing the off-diagonal terms we obtain 
$$
c_P(T_1,T_2) = \<v_1, v_2\>_\Cal E.  \tag{2.16}
$$

Now $T_k = T_{(0,v_k)}$.  We bring in the $c_k$ as follows.  
From the definition (2.2) of $c_P$ and the fact that 
$T_{(c,v)}(t) = e^{ct}T_{(0,v)}(t)$, it follows that 
$$
c_P(T_{(c_1,v_1)},T_{(c_2,v_2)}) = c_1 + \bar c_2 + 
c_P(T_{(0,v_1)},T_{(0,v_2)}).  
$$
Together with (2.16), this implies 
$$
c_P(T_{(c_1,v_1)},T_{(c_2,v_2)}) = c_1 + \bar c_2 + 
\<v_1,v_2\>_\Cal E,
$$
as required.

Once one has this formula we obtain $d_*(P)=\dim \Cal E$
by a straightforward calculation.  Indeed, the preceding arguments 
show that we may identify $\Cal U_P$ with $\Bbb C\times\Cal E$ 
in such a way that the covariance function becomes 
$$
c_P((a,v),(b,w))=a + \bar b + \<v,w\>.  
$$
Now as in \cite{1, Proposition 5.3}
one computes directly that 
$$
d_*(P)=\dim H(\Bbb C\times \Cal E,c_P) = \dim\Cal E.  
$$
\qed
\enddemo

Together with the results of \cite{4}, Theorem 
2.3 shows how to calculate the index of the minimal
\esg\ dilation of any unital CP semigroup having 
bounded generator:  

\proclaim{Corollary 2.17}
Let $L$ be a bounded linear map on $\Cal B(H)$ which 
generates a semigroup $P = \{P_t = \exp(tL): t\geq 0\}$
of normal completely positive maps satisfying 
$P_t(\bold 1) = \bold 1$
for every $t\geq 0$.  Let $\alpha = \{\alpha_t: t\geq 0\}$ 
be the minimal \esg -dilation of $P$.  Then 
$\Cal U_\alpha \neq \emptyset$ and the numerical index 
$d_*(\alpha)$ is the rank of $L$.  
\endproclaim

\demo{proof}
By Proposition 1.20, we find an 
operator $k$ and a metric 
operator space $\Cal E$ satisfying 
$\Cal E\cap \Bbb C\bold 1 = \{0\}$, and 
$$
L(x) = P_\Cal E(x) + kx + xk^*, \qquad x\in \Cal B(H).  
$$
By \cite{cpindex, Theorem 4.9} we have 
$d_*(\alpha)=d_*(P)$, while by Theorem 2.3 above
we have $d_*(p)=\dim \Cal E = \text{rank}(L)$.
\qed
\enddemo

\remark{Remarks}
It is possible, of course, for the rank of 
$L$ to be $0$; equivalently,  $\Cal E=\{0\}$ and 
hence $\Bbb C\times \Cal E \cong \Bbb C$. In this event 
$P$ is a semigroup of $*$-automorphisms of $\Cal B(H)$ 
and $\alpha = P$.  This degenerate case is discussed 
more fully  in the proof of Corollary 4.25 below.  
\endremark

\subheading{3.  Completeness of the covariance function}
It is possible for different CP semigroups $P$, $Q$ to 
have the same covariance function in the sense that 
$P$ and $Q$ have the same set of units and 
$$
c_P(S,T) = c_Q(S,T), \qquad S,T\in \Cal U_P = \Cal U_Q.  
$$
For example, let $\alpha=\{\alpha_t: t\geq 0\}$ be an 
\esg\ acting on $M=\Cal B(H)$.  For every $t>0$ let 
$\Cal E_\alpha(t)$ be the metric operator space associated
with $\alpha_t$.  Since $\alpha_t$ is an endomorphism
we have in this case 
$$
\Cal E_\alpha(t) = \{T\in \Cal B(H): \alpha_t(x)T=Tx, 
\quad x\in \Cal B(H)\}, 
$$
and the inner product in $\Cal E_\alpha(t)$ is defined 
by 
$$
\<S,T\>_{\Cal E_\alpha(t)} \bold 1 = T^*S.  
$$
Assuming that $\Cal U_\alpha\neq \emptyset$, we can 
form a closed subspace $\Cal D(t)$ of the Hilbert 
space $\Cal E_\alpha(t)$ generated by all finite 
products obtained from units as follows
$$
\Cal D(t) = [u_1(t_1)u_2(t_2)\dots u_n(t_n): 
u_k\in \Cal U_\alpha, t_k>0, t_1+\dots +t_n = t, n\geq 1\}.  
$$
$\Cal D(t)$ is itself a metric operator space, and it determines
a $*$-endomorphism $\beta_t$ of $\Cal B(H)$ by way of 
$$
\beta_t(x) = \sum_kv_kxv_k^*, \qquad x\in\Cal B(H)
$$
$\{v_1,v_2\dots\}$ denoting an orthonormal basis 
for $\Cal D(t)$.  Since $\Cal D(s+t)$ is spanned by the 
set of products $\{ST: S\in\Cal D(s), T\in\Cal D(t)\}$ it 
follows that $\beta_{s+t} = \beta_s\beta_t$. 
If we set $\beta_0=\iota_{\Cal B(H)}$, then 
$\beta = \{\beta_t: t\geq 0\}$ is a semigroup of 
normal $*$-endomorphisms of $\Cal B(H)$.  The individual 
maps $\beta_t$ are not necessarily unit preserving, but 
we do have 
$$
\beta_t\leq \alpha_t, \qquad t\geq 0,\tag{3.2}
$$
and in fact the semigroup $\beta$ is continuous.  
Moreover, the following conditions are equivalent 
for every $t>0$,
\roster
\item"{CS1}"
$H$ is spanned by $\{T\xi: T\in \Cal D(t), \xi\in H\}$
\item"{CS2}"
$\beta_t(\bold 1)=\bold 1$
\item"{CS3}"
$\Cal E_\alpha(t)=\Cal D(t)$ 
\item"{CS4}"
$\beta_t=\alpha_t$
\endroster
(see \cite{1, \S 7}), and when these conditions are satisfied for 
some $t>0$ then they are satisfied for every $t>0$.  When 
that is the case, $\alpha$ is called {\bf completely 
spatial}.  If $\alpha$ is any semigroup for which 
$\Cal U_\alpha\neq \emptyset$ then $\alpha$ is called 
{\bf spatial} and in this case we refer to its associated
semigroup $\beta$ as the {\bf standard part} of 
$\alpha$.  

Now a straightforward computation shows that $\alpha$ 
and $\beta$ have the same set of units and the same 
covariance function.  So if $\alpha$ is any spatial 
\esg\ which is {\it not} 
completely spatial, then $\alpha$, $\beta$ provide
rather extreme examples of distinct CP semigroups 
which have the same covariance function.  The existence
of such \esg s is established in \cite{10}.  
The following result asserts that this phenomenon 
cannot occur for CP semigroups which
have bounded generators.  

\proclaim{Theorem 3.3}
Let $P_1$, $P_2$ be CP semigroups with 
bounded generators $L_1$, $L_2$.  
Suppose that $P_1$ and $P_2$ have 
the same set of units and 
$$
c_{P_1}(T,T^\prime) = c_{P_2}(T,T^\prime), 
\qquad T,T^\prime\in\Cal U_{P_1}=\Cal U_{P_2}.
$$
Then $L_1=L_2$, and hence $P_1=P_2$.  
\endproclaim

\demo{poroof}
By Proposition 1.20 we can find metric operator spaces 
$\Cal E_1$, $\Cal E_2$ and operators $k_1$, $k_2\in\Cal B(H)$ 
satisfying 
$$
\Cal E_1\cap\Bbb C\bold 1=\Cal E_2\cap\Bbb C\bold 1 = \{0\}
\tag{3.4}
$$
and 
$$
L_j(x) = P_{\Cal E_j}(x) + k_jx + xk_j^*, \qquad x\in\Cal B(H)
$$
for $j=1,2$.  Theorem 2.3 asserts that the most general unit of $P_j$ 
is a semigroup of the form $T(t) = \exp(ta)$, where 
$a$ belongs to the set of operators 
$\Cal E_j + \Bbb C\bold 1 + k_j$.  Thus, the hypothesis 
$\Cal U_{P_1} = \Cal U_{P_2}$ is equivalent to the equality
of the two sets
$$
\Cal E_1 + \Bbb C\bold 1 + k_1 = \Cal E_2 + \Bbb C\bold 1+k_2.  
$$
Now if $E_1$ and $E_2$ are linear subspaces of a vector space
$V$ and $k_1$, $k_2$ are elements of $V$ such that 
$E_1 + k_1 = E_2 + k_2$ then we must have $E_1 = E_2$ and 
$k_2-k_1\in E_2$.  Taking $E_j = \Cal E_j + \Bbb C\bold 1$ 
it follows that 
$$
\Cal E_1 + \Bbb C\bold 1 = \Cal E_2 + \Bbb C\bold 1, \tag{3.6}
$$
and 
$$
k_2 = k_1 + v_2 + \lambda\bold 1, \tag{3.7}
$$
where $v_2\in \Cal E_2$ and $\lambda$ is a scalar. 

Associated with any pair of operator spaces
$\Cal E_1$, $\Cal E_2$ satisfying (3.4) and (3.6) 
there is an isomorphism of vector spaces
$\theta: \Cal E_1\to\Cal E_2$.  Indeed, since every element
$v\in \Cal E_1$ has a unique decomposition 
$$
v = v^\prime + \lambda\bold 1
$$
where $v^\prime\in \Cal E_2$ and $\lambda\in \Bbb C$, we can 
define a linear functional $f: \Cal E_1\to \Bbb C$ and a 
linear isomorphism $\theta: \Cal E_1\to \Cal E_2$ by
$$
v = \theta(v) + f(v)\bold 1.  \tag{3.8}
$$

We will show first that the linear isomorphism 
$\theta: \Cal E_1\to\Cal E_2$ defined by (3.8) is 
a unitary operator in that for any pair of 
elements $v,v^\prime\in\Cal E_1$ we have 
$\<\theta(v),\theta(v^\prime)\>_{\Cal E_2} = 
\<v,v^\prime\>_{\Cal E_1}$.  
To that end, fix
$v,v^\prime\in \Cal E_1$ and consider the units $T$, $T^\prime$ 
of $P_1$ defined by 
$$
T(t) = \exp{t(v+k_1)}, \quad T^\prime(t) = \exp{t(v^\prime + k_1)}.  
$$
Combining formula (2.2) with Theorem 2.3,
we find that for every $t>0$, 
$$
c_{P_1}(T,T^\prime) = \<v,v^\prime\>_{\Cal E_1}.  \tag{3.9}
$$
Now we must consider $T$ and $T^\prime$ relative to the coordinates
associated with $P_2$.  In order to do that, we use (3.7) and 
(3.8) to write 
$$
v+k_1 = \theta(v) + f(v)\bold 1 + k_1 = 
(\theta(v)-v_2) + (f(v)-\lambda)\bold 1 + k_2, 
$$
and similarly
$$
v^\prime + k_1 = (\theta(v^\prime)-v_2) + 
(f(v^\prime)-\lambda)\bold 1 + k_2.  
$$
Considering $T$ and $T^\prime$ as units of $P_2$, 
we have in the notation of formula (2.3.3), 
$$
T = T_{(f(v)-\lambda,\theta(v)-v_2)}, \quad
T^\prime = T_{(f(v^\prime) -\lambda,\theta(v^\prime)-v_2)}, 
$$
and corresponding to (3.9) we have 
$$
c_{P_2}(T,T^\prime) = 
f(v)-\lambda + \bar f(v^\prime) -\bar\lambda +
\<\theta(v)-v_2,\theta(v^\prime)-v_2\>_{\Cal E_2}.  \tag{3.10}
$$
Since $c_{P_1}=c_{P_2}$, we may equate 
the right sides of (2.23) and (2.24) to obtain
$$
\<v,v^\prime\>_{\Cal E_1} = 
f(v)-\lambda + \bar f(v^\prime) -\bar\lambda +
\<\theta(v)-v_2,\theta(v^\prime)-v_2\>_{\Cal E_2}.  \tag{3.11}
$$

The identity (3.11) implies that $\theta$ is unitary.  To 
see that, consider the sesquilinear form representing 
the difference 
$$
D(v,v^\prime) = \<v,v^\prime\>_{\Cal E_1} - 
\<\theta(v),\theta(v^\prime)\>_{\Cal E_2}.  
$$
We can rewrite (2.25) in the form
$$
D(v,v^\prime) = g(v) + \bar g(v^\prime), \tag{3.12}
$$
where $g: \Cal E_1\to \Bbb C$ is the function 
$g(v) = f(v) -\<\theta(v),v_2\>_{\Cal E_1}+
1/2\<v_2,v_2\>_{\Cal E_2}-\lambda $.  
For every $t>0$ we can write
$$
D(v,v^\prime) = t^{-2}D(tv,tv^\prime) = 
t^{-2}(g(tv) + \bar g(tv^\prime)), 
$$
and clearly 
$$
t^{-2}g(tv) = t^{-1}(f(v) -\<\theta(v),v_2\>) +
t^{-2}(1/2\<v_2,v_2\>_{\Cal E_2}-\lambda)
$$
tends to zero as $t\to \infty$.  Thus, 
$$
D(v,v^\prime) = 0 = g(v)+\bar g(v^\prime) \tag{3.13}
$$
for every $v,v^\prime\in\Cal E_1$.  

We claim next that (3.13) implies 
$$
\align
\lambda&=1/2\<v_2,v_2\>_{\Cal E_2}+ic, 
\qquad\text{and} \tag{3.14}\\
f(v) &= \<\theta(v),v_2\>_{\Cal E_2},\tag{3.15}
\endalign
$$
where $i=\sqrt{-1}$ and $c$ is a real number.  
Indeed, setting $v=v^\prime=0$ in the equation 
$g(v)+\bar g(v^\prime)=0$ (3.13) we obtain
$$
1/2\<v_2,v_2\>_{\Cal E_2}-\lambda +
1/2\<v_2,v_2\>_{\Cal E_2} -\bar\lambda = 0,
$$
hence (3.14).  Thus the linear functional 
$\rho(v)=f(v)-\<\theta(v),v_2\>_{\Cal E_2}$ satisfies 
$$
\rho(v) + \bar\rho(v^\prime) = g(v)+\bar g(v^\prime) = 0
$$
for all $v,v^\prime\in\Cal E_1$ and (3.15) follows 
after setting $v^\prime=0$.  

From (3.7) and (3.14) we obtain 
$$
k_2 = k_1+v_2 +1/2\<v_2,v_2\>_{\Cal E_2}\bold 1 + ic\bold 1,
$$
so for all $x\in\Cal B(H)$ we have 
$$
k_2x+xk_2^* = k_1x+xk_1^* + 
v_2x+xv_2^*+\<v_2,v_2\>_{\Cal E_2}x.  
\tag{3.16}
$$
Since $L_2(x) = P_{\Cal E_2}(x) + k_2x+xk_2^*$ it follows that 
$$
L_2(x)=P_{\Cal E_2}(x)+v_2x+xv_2^*+
\<v_2,v_2\>_{\Cal E_2}x +k_1x+xk_1^*.
\tag{3.17}
$$
We will show now that the right side of (3.17) is $L_1(x)$; 
equivalently, we claim that 
$$
P_{\Cal E_2}(x) +v_2x+xv_2^* +
\<v_2,v_2\>_{\Cal E_2}x = P_{\Cal E_1}(x). \tag{3.18}
$$
To see that, let $u_1, u_2,\dots$ be an orthonormal basis 
for $\Cal E_1$.  Then $\theta(u_1),\theta(u_2)\dots$ is an 
orthonormal basis for $\Cal E_2$ and if we set 
$\mu_k=\<v_2,\theta(u_k)\>_{\Cal E_2}$ then the sequence 
$(\mu_1,\mu_2,\dots)$ belongs to $\ell^2$ and 
$v_2=\sum_k\mu_k\theta(u_k)$.  Thus we have 
$$
v_2x+xv_2^*+\<v_2,v_2\>_{\Cal E_2}x = 
\sum_k(\mu_k\theta(u_k)x+x\bar\mu_k\theta(u_k)^* + |\mu_k|^2x),
$$
while 
$$
P_{\Cal E_2}(x) = \sum_k\theta(u_k)x\theta(u_k)^*,
$$
so that the left side of (3.18) can be written 
$$
\sum_k(\theta(u_k)+\bar\mu_k\bold 1)x(\theta(u_k)+\bar\mu_k\bold 1)^*.
$$
Noting that for each $k$, 
$$
\theta(u_k)+\bar\mu_k\bold 1 = 
\theta(u_k)+\<\theta(u_k),v_2\>_{\Cal E_2}\bold 1 =
\theta(u_k)+f(u_k)\bold 1 = u_k
$$
by the definition (3.8) of $\theta$ and $f$, we find that the last 
expression reduces to
$$
\sum_k u_kx u_k^* = P_{\Cal E_1}(x),
$$
as asserted.  That completes the proof of Theorem 3.3\qed
\enddemo

\remark{Remark 3.18}
The proof of Theorem 3.3 gives somewhat more information 
than is contained in its statement.  For example, suppose that 
$\Cal E_1$, $\Cal E_2$ are two metric operator spaces 
satisfying $\Cal E_1\cap\Bbb C\bold 1=\Cal E_2\cap\bold 1=\{0\}$, 
and let $k_1$, $k_2$ be two operators such that the corresponding 
generators are the same:
$$
P_{\Cal E_1}(x)+k_1x+xk_1^* = P_{\Cal E_2}+k_2x+xk_2^*, 
\qquad x\in\Cal B(H).  
$$
Then the proof of Theorem 3.3 implies that 
$$
\Cal E_1+\Bbb C\bold 1 = \Cal E_2+\Bbb C\bold 1.  \tag{3.19}
$$
(3.19) allows one to define a linear isomorphism 
$\theta: \Cal E_1\to\Cal E_2$ and a linear functional $f$ 
on $\Cal E_1$ by
$$
v=\theta(v)+f(v)\bold 1,\qquad v\in\Cal E_1. \tag{3.20}
$$
The same argument shows that $\theta$ is a unitary operator 
and that the unique element $v_2\in\Cal E_2$ defined by 
$$
f(v)=\<\theta(v),v_2\>_{\Cal E_2},\qquad v\in\Cal E_1 \tag{3.21}
$$
satisfies
$$
k_2 = k_1+v_2+(1/2\<v_2,v_2\>_{\Cal E_2} + ic)\bold 1 \tag{3.22}
$$
where $c$ is a real constant.  

Conversely, if we start with two pairs $(\Cal E_1,k_1)$,  
$(\Cal E_2,k_2)$ satisfying (3.19)--(3.22) (so that 
the linear isomorphism $\theta$ defined by (3.20) is
a unitary operator) along with 
$\Cal E_1\cap\Bbb C\bold 1=\Cal E_2\cap\Bbb C\bold1=\{0\}$, 
then both $(\Cal E_1,k_1)$ and $(\Cal E_2,k_2)$ 
give rise to the same generator
$$
L(x) = P_{\Cal E_1}(x)+k_1x+xk_1^* = P_{\Cal E_2}(x)+k_2x+xk_2^*.
$$
\endremark

\subheading{4. Minimal dilations}

Every unital CP semigroup $P = \{P_t: t\geq 0\}$ acting on 
$\Cal B(H)$ has a minimal
dilation to an \esg\ $\alpha = \{\alpha_t: t\geq 0\}$
acting on $\Cal B(K)$, where $K$ is a Hilbert space 
which contains $H$ as a subspace \cite{2,5,6}.  
Recall that an \esg\ $\alpha$ 
is said to be completely spatial if there are ``enough" 
units in the sense that the equivalent conditions 
CS1--CS4 of the preceding section are satisfied.  
The completely spatial \esg s constitute the best-understood 
class.  All of the basic examples are of this type,
and they are classified up to cocycle 
conjugacy by their numerical index $d_*(\alpha)$ 
\cite{1, Corollary of Proposition 7.2}.  
Thus, if one knows that an \esg\ $\alpha$ is completely
spatial and has numerical index $n = 1, 2, \dots, \infty$, 
then $\alpha$ 
must be conjugate to a cocycle perturbation of the $CAR/CCR$ flow of
index $n$.  

The purpose of this section is to show that if a CP semigroup
has a {\it bounded} generator then its minimal dilation is 
completely spatial (Theorem 4.8 below).  The proof of 
Theorem 4.8 is based on Theorem 3.3 and the following 
result.  

\proclaim{Theorem 4.1}
Let $P=\{P_t: t\geq 0\}$ be a CP semigroup acting on
$\Cal B(H)$ which has a bounded generator.  Let 
$\{Q_t: t>0\}$ be a family of normal completely 
positive maps on $\Cal B(H)$ satisfying the two 
conditions
$$
\align
Q_t &\leq P_t, \tag{4.1.1}\\
Q_{s+t} &=Q_sQ_t,\tag{4.1.2}
\endalign
$$
for all $s,t>0$, and which is not the trivial family
$Q_t=0$, $t>0$.  Let $Q_0$ be the identity map of $\Cal B(H)$.  
Then $\{Q_t: t\geq 0\}$ is also a CP semigroup having 
bounded generator.  
\endproclaim

Our proof of this result requires the following estimate.  

\proclaim{Lemma 4.2}
Let $P$ be a normal completely positive linear map on 
$M = \Cal B(H)$ and let $\sigma_P$ be its symbol.  
Then we have
$$
\inf_{\lambda>0}\|P - \lambda\iota_M\| \leq 6\|P\|^{1/2}
\sup_{\|x\|\leq 1}\|\sigma_P(dx^*\,dx)\|^{1/2}, 
$$
where $\iota_M$ denotes the identity map of $M$.  
\endproclaim

\remark{Remark}
In proving this estimate we will make use of the following 
bit of lore.  Let $N\subseteq\Cal B(K)$ be an amenable 
von Neumann algebra and let $T$ be an operator
on $K$.  Then there is an operator $T^\prime$ in the 
commutant of $N$ such that 
$$
\|T-T^\prime\| = \sup\{\|Tx - xT\|: x\in N, \|x\|\leq 1\}.  \tag{4.3}
$$
Indeed, the operator $T^\prime$ may be obtained by a 
familiar averaging process, in which one uses 
an invariant mean on 
a suitable subgroup $G$ of the unitary group in $N$
to average quantities of the form $uTu^*$, $u\in G$, 
after noting 
that for every such unitary operator $u$, 
$\|uTu^* - T\|$ is dominated by the right side 
of (4.3).  
\endremark
\demo{proof of Lemma 4.2}
Because of Stinespring's theorem there is a Hilbert space 
$K$, a normal representation $\pi: M\to \Cal B(K)$ and an 
operator $V\in \Cal B(H,K)$ such that $P(x) = V^*\pi(x)V$, 
$x\in M$.  As in the proof of Lemma 1.19, the symbol of 
$P$ is related to $V$ and $\pi$ by way of 
$$
(\pi(x)V - Vx)^*(\pi(x)V - Vx) = 
\sigma_P(dx^*\,dx), 
$$
and hence
$$
\|Vx - \pi(x)V\|^2 = \|\sigma_P(dx^*\,dx)\|.  
$$
Setting $\|\sigma_P\| = \sup_{\|x\|\leq 1}\|\sigma_P(dx^*\,dx)\|$, 
we obtain
$$
\sup_{\|x\|\leq 1}\|Vx - \pi(x)V\|\leq \|\sigma_P\|^{1/2}.  \tag{4.4}
$$
Now consider the von Neumann algebra $N\subseteq\Cal B(H\oplus K)$
and the operator $T\in \Cal B(H\oplus K)$ defined by 
$$
N =\{ 
\pmatrix
x & 0\\
0 & \pi(x)
\endpmatrix 
: x\in M\}, \qquad
T = \pmatrix
0 & V^*\\
V & 0
\endpmatrix 
.  
$$
One finds that 
$$
T
\pmatrix
x & 0\\
0 & \pi(x)
\endpmatrix
-
\pmatrix
x & 0\\
0 & \pi(x)
\endpmatrix
T = 
\pmatrix
0 & -(Vx - \pi(x)V)^*\\
Vx - \pi(x)V & 0
\endpmatrix
.
$$
The norm of the operator on the right is 
$\|Vx - \pi(x)V\|$, hence
$$
\align
&\sup\{\|Ty - yT\|: y\in N, \|y\|\leq 1\} = \\
&\sup\{\|Vx - \pi(x)V\|: x\in M, \|x\|\leq 1\} \leq \|\sigma_P\|^{1/2}.  
\endalign
$$
Using (4.3) we find an operator $T^\prime\in N^\prime$ satisfying
$\|T - T^\prime\| \leq \|\sigma_P\|^{1/2}$.  A straightforward
computation shows that operators in the commutant of $N$ must
have the form
$$
T^\prime = 
\pmatrix
A & Y^* \\
X & B
\endpmatrix
$$
where $A$ is a scalar multiple of the identity of $\Cal B(H)$, 
$B$ belongs to the commutant of $\pi(M)$, and $X$ and $Y$ 
are intertwining operators, $Xx = \pi(x)X$, $Yx = \pi(x)Y$, 
$x\in M$.  It follows that there is such an $X$ for which 
$$
\|V - X\|\leq \|T - T^\prime\| \leq \|\sigma_P\|^{1/2}. \tag{4.5}
$$  
Since $X^*X$ commutes with $M=\Cal B(H)$ we must have 
$X^*X = \lambda \bold 1_H$ for some scalar $\lambda\geq 0$, 
and hence 
$$
\|P(x) - \lambda\cdot x\| = \|V^*\pi(x)V - X^*\pi(x)X\| \leq 
2\|V - X\|\cdot\|x\|\max{(\|V\|,\|X\|)}.  
$$ 
Note that $\max{(\|V\|,\|X\|)}\leq 3\|P\|^{1/2}$.  Indeed, 
since $V^*V = P(\bold 1)$ we have $\|V\|\leq \|P\|^{1/2}$, and
by (1.4) we can estimate $\|X\|$ by way of 
$$
\|X\|\leq \|V\| + \|V - X\| \leq \|P\|^{1/2} + \|\sigma_P\|^{1/2} \leq
3\|P\|^{1/2}.  
$$ 
Using (4.5) we arrive at the desired inequality
$\|P - \lambda\cdot\iota_M\| \leq 6\|P\|^{1/2}\|\sigma_P\|^{1/2}$
\qed
\enddemo

\demo{proof of Theorem 4.1}
Let $\{Q_t: t>0\}$ satisfy (4.1.1) and (4.1.2).  We will 
show that 
$$
\lim_{t\to 0+}\|Q_t - \iota\| = 0, \tag{4.6}
$$
$\iota$ denoting the identity map $\iota(x)=x$, $x\in\Cal B(H)$.  
Theorem 4.1 follows immediately since under these conditions 
the semigroup $\{Q_t: t\geq 0\}$ becomes a 
continuous semigroup of elements in the Banach algebra 
of all normal linear mappings $L$ on $\Cal B(H)$
with the uniform norm
$$
\|L\| = \sup_{\|x\|\leq 1}\|L(x)\|.  
$$

In order to prove (4.6), we claim first 
that there is a family $\lambda_t$, $t>0$ 
of nonnegative numbers such that $\|Q_t-\lambda_t\iota\|\to 0$
as $t\to 0+$.  Indeed, since $P_t-Q_t$ is completely positive
for every $t>0$ we 
have $0\leq \sigma_{Q_t}(dx\,dx^*) \leq \sigma_{P_t}(dx\,dx^*)$ 
for every $x\in M$; hence 
$$
\|\sigma_{Q_t}(dx\,dx^*)\| \leq \|\sigma_{P_t}(dx\,dx^*)\| =
\|\sigma_{P_t-\iota}(dx\,dx^*)\| \leq 4\|P_t-\iota\|\,\|x\|^2.  
$$
Using Lemma 4.2 together with the latter inequality we find that 
$$
\align
\inf_{\lambda>0}\|Q_t - \lambda\iota\|\leq 
6&\|Q_t\|^{1/2}\sup_{\|x\|\leq 1}\|\sigma_{Q_t}(dx\,dx^*)\|^{1/2} \leq \\
12&\|Q_t\|^{1/2}\|P_t-\iota\|^{1/2} \leq 
12\|P_t\|^{1/2}\|P_t-\iota\|^{1/2}
\endalign
$$
and the claim follows because $\|P_t-\iota\|$ tends to 
$0$ as $t\to 0+$.  

It remains to prove that $\lambda_t\to 1$ as $t\to 0+$.  
To that end, we claim 
$$
\lim_{t\to 0+}\|Q_t\| = 1.  \tag{4.7}
$$
Indeed, since 
$\|Q_t\| \leq \|P_t\|$ and $\|P_t - \iota\|\to 0$ as 
$t\to 0+$, we have 
$$
\limsup_{t\to 0+}\|Q_t\|\leq \limsup_{t\to0+}\|P_t\|=\|\iota\|= 1.  
$$
So if (4.7) fails then we must have 
$$
\liminf_{t\to 0+}\|Q_t\|  < 1, 
$$
and in that event we can pick $r<1$ such that 
$\liminf_{s\to 0+}\|Q_s\|\ < r$.  
Let $R>1$ be close enough to $1$ so that 
$rR < 1$.  Then for sufficiently small $t$ we have 
both $\|Q_t\|\leq R$ and $\inf_{0<s<t}\|Q_s\|\leq r$.  
For such a $t$ we can find $0<s<t$ such that 
$\|Q_s\|\leq r$, and hence
$$
\|Q_t\| \leq \|Q_s\|\cdot\|Q_{t-s}\|\leq rR.  
$$
Thus $\|Q_t\|\leq \|Q_{t/n}\|^n\leq (rR)^n$ for every
$n= 1,2,\dots$ and hence $\|Q_t\| = 0$ for all sufficiently
small $t$.  Because of the semigroup property 
it follows that $\|Q_t\|$ is identically zero, 
contradicting our 
hypothesis on $Q$ and proving (4.7).   

To see that $\lambda_t\to 1$ as $t\to 0+$ write
$$
\lambda_t = \|(\lambda_t\iota-Q_t) + Q_t\|, 
$$
and use $\lim_{t\to 0+}\|Q_t - \lambda_t\iota\|=0$ together 
with (4.7)\qed
\enddemo

Following is the principal result of this section.  

\proclaim{Theorem 4.8}
Let $P=\{P_t: t\geq 0\}$ be a unital CP semigroup having
a bounded generator.  Then the minimal dilation of $P$ 
to an \esg\ is completely spatial.  
\endproclaim
\demo{proof}
Let $\alpha=\{\alpha_t: t\geq 0\}$ be an \esg\ acting 
on $\Cal B(K)$, $K\supseteq H$ which is the minimal dilation
of $P$.  Let $\beta = \{\beta_t: t\geq0\}$ be the standard
part of $\alpha$ as defined in section 3.  Letting 
$\{\Cal D(t): t>0\}$ be the family of metric operator 
spaces of (3.1), we have to show that 
$$
K = [A\xi: A\in\Cal D(t), \xi\in K],    
$$
or equivalently, that $\beta_t(\bold 1_K)=\bold 1_K$ for 
every $t$.  

Now in general, since $\beta$ is a semigroup of endomorphisms
of $\Cal B(K)$, the projections $\beta_t(\bold 1_K)$ decrease
as $t$ increases, and the limit projection 
$$
e_\infty = \lim_{t\to\infty}\beta_t(\bold 1_K) 
\tag{4.9}
$$
is fixed under the action of $\beta$.  The compression 
$\beta^\infty$ of the semigroup $\beta$ to the corner 
$e_\infty\Cal B(K)e_\infty$ can be 
viewed as a semigroup of {\it unit-preserving}
endomorphisms of $\Cal B(e_\infty K)$.  It is not quite
obvious that $\beta^\infty$ is an \esg\ since we have 
not proved that $\beta$ is continuous.  While it is 
possible to establish that directly, 
we will not have to do so since the following 
result implies that $\beta^\infty$ is actually a compression
of $\alpha$.  

\proclaim{Lemma 4.10}
The projection $e_\infty$ satisfies 
$\alpha_t(e_\infty)\geq e_\infty$ for every $t\geq 0$, and 
the compression of $\alpha_t$ to $e_\infty\Cal B(K)e_\infty$ 
is $\beta^\infty_t$ for every $t\geq 0$.  In particular, 
$\beta^\infty$ defines an \esg\ acting on $\Cal B(e_\infty K)$.    
\endproclaim
\demo{proof}
We claim first that for every $x\in\Cal B(K)$ we have 
$$
\alpha_t(x)\beta_t(\bold 1_K) = 
\beta_t(x)=\beta_t(\bold 1_K)\alpha_t(x).  \tag{4.11}
$$
Indeed, if we let $u_1, u_2,\dots$ be an orthonormal 
basis for $\Cal D(t)$ then we have 
$$
\beta_t(x) = \sum_j u_jxu_j^*,\qquad x\in\Cal B(K).  
$$
Since $\Cal D(t)\subseteq\Cal E_\alpha(t)$ it follows 
that for every $j$,
$$
\alpha_t(x)u_ju_j^* = u_jxu_j^* = u_ju_j^*\alpha_t(x),
$$
and (4.11) follows by summing on $j$.  
Taking $x=e_\infty$ in (4.11) and using 
$e_\infty = \beta_t(e_\infty)e_\infty$, we obtain
$$
\alpha_t(e_\infty)e_\infty =
\alpha_t(e_\infty)\beta_t(e_\infty)e_\infty = 
\beta_t(e_\infty)e_\infty = e_\infty,
$$
hence $\alpha_t(e_\infty)\geq e_\infty$.  

Now choose an operator $x\in e_\infty\Cal B(K)e_\infty$.  
We have
$$
\align
e_\infty\alpha_t(x)e_\infty &= 
e_\infty\beta_t(\bold 1_K)\alpha_t(x)e_\infty =
e_\infty\beta_t(x)e_\infty \\
&=\beta_t(e_\infty)\beta_t(x)\beta_t(e_\infty) = 
\beta_t(e_\infty x e_\infty) = \beta_t(x),
\endalign
$$
as asserted \qed
\enddemo

Let $p_0\in\Cal B(K)$ be the projection onto $H$.  Then of 
course $\alpha_t(p_0)\geq p_0$, and after identifying 
$\Cal B(H)$ with $p_0\Cal B(K)p_0$ we have 
$$
P_t(x) = p_0\alpha_t(x)p_0, \qquad t\geq 0, x\in\Cal B(H).  
\tag{4.12}
$$
We will prove that 
$$
p_0\leq e_\infty = \lim_{t\to\infty}\beta_t(\bold 1_K).  \tag{4.13}
$$
Granting (4.13) for a moment, then Lemma 4.10 asserts that 
the compression of $\alpha$ to $e_\infty\Cal B(K)e_\infty$ is
an \esg\ which is obviously an intermediate dilation of 
$P$.  By the {\it definition} of minimal dilation given
in \cite{2} we may conclude that $e_\infty=\bold 1_K$.  
In particular, $\beta_t(\bold 1_K)=\bold 1_K$ for every $t$, 
proving that $\alpha$ is completely spatial.  

Thus the proof of Theorem 4.8 is reduced to establishing 
(4.13).  This will be done indirectly, by considering the 
family of maps $Q=\{Q_t: t\geq 0\}$ 
on $\Cal B(H)$ obtained by compressing the semigroup 
$\beta$:
$$
Q_t(x) = p_0\beta_t(x)p_0, \qquad x\in\Cal B(H), t\geq 0.
\tag{4.14}
$$
We will show that $Q$ is a CP semigroup with bounded
generator, which has the same set of units and the same 
covariance function as $P$.  Theorem 3.3 will then imply 
that $Q_t=P_t$ for all $t$, and in particular 
$Q_t(p_0)=P_t(P_0) =p_0$.  From the definition (4.14) 
of $Q$ we conclude that 
$$
p_0\leq \beta_t(p_0) \leq \beta_t(\bold 1_K)
$$
for every $t\geq 0$, and (4.13) will follow.  

We show first that the family $Q$ of (4.14) is a semigroup.  
To that end, we claim that for every $t\geq 0$, 
$$
p_0\beta_t(p_0) = p_0\beta_t(\bold 1_K),\quad\text{and }
\beta_t(p_0)p_0=\beta_t(\bold 1_K)p_0.  \tag{4.15}
$$
Indeed, since 
$\beta_t(x)=\alpha_t(x)\beta_t(\bold 1_K)$ for every 
$x\in\Cal B(K)$ we have 
$$
p_0\beta_t(\bold 1_K-p_0)=
p_0\alpha_t(\bold 1_K-p_0)\beta_t(\bold 1_K)=
p_0(\bold 1_K-\alpha_t(p_0))\beta_t(\bold 1_K)
$$
and the right side is zero because $\alpha_t(p_0)\geq p_0$.
Similarly, $\beta_t(\bold 1_K-p_0)p_0 = 0$.  
Thus for $s,t>0$ and $x\in\Cal B(H)$ we have 
$$
Q_sQ_t(x) = p_0\beta_s(p_0\beta_t(x)p_0)p_0 =
p_0\beta_s(p_0)\beta_{s+t}(x)\beta_s(p_0)p_0.  
$$
Because of (4.15) we have 
$p_0\beta_s(p_0)=p_0\beta_s(\bold 1_K)$ and 
$\beta_s(p_0)p_0=\beta_s(\bold 1_K)p_0$, hence 
$$
Q_sQ_t(x)=
p_0\beta_s(\bold 1_K)\beta_{s+t}(x)\beta_s(\bold 1_K)p_0.  
$$
Finally, since $\beta_{s+t}(\bold 1_K)\leq \beta_s(\bold 1_K)$, 
we have 
$$
\beta_s(\bold 1_K)\beta_{s+t}(x)\beta_s(\bold 1_K) =\beta_{s+t}(x), 
$$
and the formula $Q_sQ_t(x)=Q_{s+t}(x)$ follows.  

Now since $\beta_t\leq \alpha_t$ for every $t\geq 0$ we have 
$Q_t\leq P_t$.  Thus the hypotheses (4.1.1) and (4.1.2) of 
Theorem 4.1 are satisfied.  Notice that $Q_t$ cannot vanish 
for every $t>0$.  Indeed, Corollary 2.17 implies that $\alpha$ 
must have units.  Choose $u\in\Cal U_\alpha$, and let $T$ be 
the corresponding unit of $P$ defined by
$$
T(t)^* = u(t)^*\restriction_H, \qquad t\geq 0.  
$$
There is a real constant $k$ such that 
$\<u(t),u(t)\>_{\Cal E_\alpha(t)} = e^{kt}$.  Since 
$u(t)\in\Cal D(t)$ we have 
$$
u(t)u(t)^*\leq e^{kt}\beta_t(\bold 1_K),
$$
and hence by (4.15) 
$$
T(t)T(t)^*=p_0u(t)u(t)^*p_0\leq e^{kt}p_0\beta_t(\bold 1_K)p_0 
=e^{kt}p_0\beta_t(p_0)p_0 = Q_t(p_0).  
$$
Since $T(t)T(t)^*$ tends weakly to $p_0=\bold 1_H$ as $t\to0+$, 
so does $Q_t(p_0)$.  

From Theorem 4.1 we conclude that the semigroup $Q$ has a bounded 
generator.  Finally, we claim that $Q_t=P_t$ for all $t\geq 0$.  
According to Theorem 3.3, it is enough to show that 
$\Cal U_P=\Cal U_Q$ and that $P$ and $Q$ have the same covariance
function.  For that, we require a general observation: 

\proclaim{Lemma 4.16}
Let $P$, $Q$ be two CP semigroups acting on $\Cal B(H)$, and 
suppose that $Q_t\leq P_t$ for every $t\geq 0$.  Then 
$\Cal U_Q\subseteq \Cal U_P$, and for every finite set 
$T_1,T_2,\dots,T_n\in\Cal U_Q$ and 
$\lambda_1,\lambda_2,\dots,\lambda_n\in\Bbb C$ we have 
$$
\sum_{i,j}\lambda_i\bar\lambda_jc_P(T_i,T_j) \leq 
\sum_{i,j}\lambda_i\bar\lambda_jc_Q(T_i,T_j), 
$$
where $c_P$ and $c_Q$ are the covariance functions of 
$P$ and $Q$.  
\endproclaim
\remark{Remarks}
Let $X$ be a set and $f: X\times X\to\Bbb C$ a function.  
Recall that $f$ is {\it positive definite} if for every 
$x_1,x_2,\dots,x_n\in X$ and 
$\lambda_1,\lambda_2,\dots,\lambda_n\in\Bbb C$ we have 
$$
\sum_{i,j}\lambda_i\bar\lambda_jf(x_i,x_j)\geq 0.  
$$
Given two functions $f,g: X\times X\to\Bbb C$, we will 
write $f\lesssim g$ or $g\gtrsim f$ if $g-f$ is positive definite.  
We will make use of the following 
elementary facts about positive 
definite functions.  
$$
f\gtrsim 0,\quad g\gtrsim 0 \implies fg\gtrsim 0, \tag{4.17}
$$
where $fg$ denotes the pointwise product $fg(x,y)=f(x,y)g(x,y)$.  
(4.13), together with transitivity of the relation $\lesssim$,
implies that for any four complex-valued functions 
$f_1,f_2,g_2,g_2$ on $X\times X$ we have 
$$
0\lesssim f_i\lesssim g_i, i=1,2\implies f_1f_2\lesssim g_1g_2.  
\tag{4.18}
$$
\endremark
\demo{proof of Lemma 4.16}
Choose $T\in\Cal U_Q$.  Then there is a real constant 
$k$ such that the semigroup $R_t(x)=T(t)xT(t)^*$ satisfies
$$
R_t\leq e^{kt}Q_t\leq e^{kt}P_t, \qquad t\geq 0,
$$
hence $T\in\Cal U_P$.  

Fix $S,T\in\Cal U_Q$, and for every $t>0$ define functions
$f,g$ by
$$
\align
f(S,T;t) &= \<S(t),T(t)\>_{\Cal E_Q(t)}, \\
g(S,T;t) &= \<S(t),T(t)\>_{\Cal E_P(t)}.  
\endalign
$$
Notice first that 
$$
0\lesssim f(\cdot,\cdot;t)\lesssim g(\cdot,\cdot;t)  \tag{4.19}
$$
as functions on $\Cal U_Q\times \Cal U_Q$.  Indeed, if 
$T_1,T_2,\dots,T_n\in\Cal U_Q$, 
$\lambda_1,\lambda_2,\dots,\lambda_n\in\Bbb C$ and we set 
$A(t)=\sum_k\lambda_kT_k(t)$, then $A(t)$ belongs to 
$\Cal E_Q(t)$ and hence the mapping 
$x\mapsto A(t)xA(t)^*$ is dominated by 
$$
\<A(t),A(t)\>_{\Cal E_Q(t)} Q_t\leq \<A(t),A(t)\>_{\Cal E_Q(t)}P_t.  
$$
It follows that 
$$
0\leq\<A(t),A(t)\>_{\Cal E_P(t)}\leq \<A(t),A(t)\>_{\Cal E_Q(t)}
$$
and (4.19) follows after expanding the inner products in 
the obvious way.  

For every partition $\Cal P=\{0=t_0<t_1,\dots<t_n=t\}$ of the 
interval $[0,t]$, set 
$$
\align
f_{\Cal P}(S,T;t) &=
\prod_{k-1}^n\<S(t_k-t_{k-1}),T(t_k-t_{k-1})\>_{\Cal E_P(t)}, \\
g_{\Cal P}(S,T;t) &=
\prod_{k-1}^n\<S(t_k-t_{k-1}),T(t_k-t_{k-1})\>_{\Cal E_Q(t)}.  
\endalign
$$
By (4.17) and (4.18), we have 
$$
0\lesssim f_{\Cal P}(\cdot,\cdot;t)\lesssim g_{\Cal P}(\cdot,\cdot; t)
$$
for every partition $\Cal P$ of $[0,t]$.  After taking the limit
on $\Cal P$ we obtain
$$
e^{tc_P}\lesssim e^{tc_Q}, 
$$
for every $t>0$.  It follows that 
$e^{tc_P}-1\lesssim e^{tc_Q}-1$ for every $t>0$ and hence
$$
c_P = \lim_{t\to 0+}\frac{e^{tc_P} - 1}{t} \lesssim 
\lim_{t\to 0+}\frac{e^{tc_Q} - 1}{t} = c_Q,
$$
as required\qed
\enddemo

We claim now that 
$$
\align
\Cal U_P&\subseteq \Cal U_Q, \qquad\qquad \text{and}\tag{4.19}\\
c_Q(T,T)&\leq c_P(T,T), \qquad \text{for every }T\in\Cal U_P.\tag{4.20}
\endalign
$$
To see that, choose any unit $T\in\Cal U_P$.  By \cite{4, Theorem 3.6}
there is a unit $u\in\Cal U_\alpha$ such that 
$$
T(t)^*= u(t)^*\restriction_H, \qquad t>0; 
$$
moreover, $c_P(T,T)=c_\alpha(u,u)$ because of the minimality of $\alpha$.  
Since $u(t)\in\Cal D(t)$ and $\<u(t),u(t)\>_{\Cal D(t)}=e^{tc_\alpha(u,u)}$,
the map 
$$
x\in\Cal B(H)\mapsto e^{tc_\alpha(u,u)}Q_t(x)-T(t)xT(t)^*
$$
is completely positive.  (4.19) follows.  We may also conclude 
from this argument that 
$$
\<T(t),T(t)\>_{\Cal E_Q(t)}\leq e^{tc_P(T,T)}.  
$$
Hence for every partition $\Cal P=\{0=t_0<t_1<\dots<t_n=t\}$
of the interval $[0,t]$, 
$$
\prod_{k=1}^n\<T(t_k-t_{k-1}),T(t_k-t_{k-1})\>_{\Cal E_Q(t)}\leq
e^{tc_P(T,T)},
$$
and after passing to the limit on $\Cal P$ we obtain
$$
e^{tc_Q(T,T)}\leq e^{tc_P(T,T)}
$$
for every $t>0$, from which (4.20) is immediate.  

Together with Lemma 4.12, (4.19) implies that $\Cal U_P=\Cal U_Q$.  
We claim now that $c_P=c_Q$.  To see that, fix $T_1$, 
$T_2\in\Cal U_P$, and consider the $2\times 2$ matrix 
$A = (a_{ij})$ defined by 
$$
a_{ij} = c_Q(T_i,T_j)-c_P(T_i,T_j).  
$$
Lemma 4.16 implies that $A\geq 0$; while (4.20) implies that 
both diagonal terms of $A$ are nonpositive, so that the trace 
of $A$ is nonpositive.  It follows that $A=0$.  In particular,
$$
c_Q(T_1,T_2)-c_P(T_1,t_2) = a_{ij} = 0.  
$$
We may now apply Theorem 3.3 to obtain $P=Q$.  As we have 
already pointed out, one may deduce from this the required 
inequality (4.13).  That completes the proof of Theorem 4.8
\qed
\enddemo

\proclaim{Corollary 4.21}
Let $P=\{P_t: t\geq 0\}$ be a unital CP semigroup acting on 
$\Cal B(H)$ whose generator is bounded and which is not 
a semigroup of $*$-automorphisms.  Then the minimal dilation 
of $P$ is a cocycle perturbation of a $CAR/CCR$ flow of 
positive index $r=1,2,\dots,\infty$.  The index $r$ is the rank 
of the generator of $P$.  
\endproclaim
\demo{proof}
By Theorem 4.8, the minimal dilation of $P$ is a completely 
spatial \esg.  The classification results of 
\cite{1, Corollary of Proposition 7.2} imply that
this \esg\ is cocycle conjugate to a $CAR/CCR$ flow.  
Its index is the rank of the generator of $P$ by 
Corollary 2.17.  

Thus we only have to check that the generator cannot have 
rank zero.  Let $L$ be the generator of $P$ and suppose that 
$L$ has rank zero.  Then the metric operator space 
associated with $L$ is $\{0\}$, and $L$ must have the form 
$$
L(x) = kx + xk^*
$$
for some $k\in M_n(\Bbb C)$ (see Definition 1.23).  
Since $P$ is unital we have $L(\bold 1) =0$,
hence $k+k^*=0$.  It follows that $k=ih$ for some self adjoint 
matrix $h$ and $L(x) = (ih)x-x(ih)$.  Thus
$$
P_t(x) = e^{ith}xe^{-ith}
$$
is a semigroup of $*$-automorphisms, contrary to hypothesis \qed
\enddemo

\remark{Remarks}
In particular, Corollary 4.21 leads to a description of 
the minimal dilations of all unital CP semigroups which 
act on a matrix algebra $M_n(\Bbb C)$, $n=2,3,\dots$.  If the 
semigroup is nontrivial then its minimal dilation $\alpha$ is 
a cocycle perturbation of a $CAR/CCR$ flow of finite positive 
index $d_*(\alpha)$.  Considering the relation between 
generators and metric operator spaces (Proposition 1.20), 
we find that for fixed $n$ the possible values of 
$d_*(\alpha)$ are $1,2,\dots,n^2-1$.  
\endremark

\end
\Refs
%
\ref\no 1\by Arveson, W.\paper Continuous analogues of Fock space
\jour Memoirs Amer. Math. Soc.\vol 80 no. 3\yr 1989
\endref

\ref\no 2\bysame \paper Minimal \esg s 
\paperinfo to appear
\endref

\ref\no 3\bysame \paper Noncommutative flows I: dynamical 
invariants \paperinfo preprint (November 1995)
\endref

\ref\no 4\bysame \paper The index of a quantum
dynamical semigroup \paperinfo preprint (July 1996)
\endref

\ref\no 5\by Bhat, B. V. R. \paper Minimal dilations of 
quantum dynamical semigroups to semigroups of endomorphisms of 
\cstar s \jour Trans. A.M.S. \toappear
\endref

\ref\no 6\bysame \paper On minimality of Evans-Hudson flows
\jour (preprint) 
\endref

\ref\no 7\by Christensen, E. and Evans, D. \paper 
Cohomology of operator algebras and quantum dynamical 
semigroups \jour J. London Math. Soc. \vol 20 
\yr 1979 \pages 358--368
\endref


\ref\no 8 \by Evans, D. and Lewis, J. T. \paper Dilations
of irreversible evolutions in algebraic quantum theory
\jour Comm. Dubl. Inst. Adv. Studies, Ser A\vol 24\yr 1977
\endref

\ref\no 9\by Powers, R. T. 
\paper A non-spatial continuous semigroup os $*$-endomorphisms 
of $\Cal B(H)$\jour Publ. RIMS (Kyoto University)\vol23\yr 1987
\pages 1053--1069
\endref

\ref\no 10\bysame \paper New examples of continuous spatial 
semigroups of endomorphisms of $\Cal B(H)$ 
\paperinfo Jour. Funct. Anal. (to appear)
\endref

\ref\no 11\bysame \paperinfo Manuscript in preparation
\endref

\ref\no 12\book Functional Analysis vol I\by Reed, M. and 
Simon, B. 
\publ Academic Press \yr 1990
\endref


\endRefs
